\documentclass[12pt]{amsart}
\usepackage[all, cmtip]{xy}

\usepackage{times,amsfonts,amsmath,amstext,amsbsy,amssymb,
amsopn,amsthm,upref,eucal, mathrsfs}

\newcommand \tr {{\mathrm{tr}}}

\newcommand \Prob {{\mathbb P}}

\newcommand \ee {{\mathbb E}}
\newcommand \Var {{\mathrm {Var}}}

\newcommand \Conf {{\mathrm {Conf}}}

\newcommand\scrI{{\mathscr I}}
\newcommand \A {{\mathscr A}}

\newtheorem{theorem}{Theorem}[section]

\newtheorem{lemma}[theorem]{Lemma}
\newtheorem{assumption}{Assumption}
\newtheorem{corollary}[theorem]{Corollary}
\newtheorem{proposition}[theorem]{Proposition}

\begin{document}

\title[Quasi-Symmetries of Determinantal Point Processes] {Quasi-Symmetries of Determinantal Point Processes.}

\begin{abstract}
The main result of this paper is that determinantal point processes on ${\mathbb R}$
corresponding to projection operators
with integrable kernels are quasi-invariant, in the continuous case,
under the group of diffeomorphisms with compact support  (Theorem \ref{main-cont});
in the discrete case, under the group of all finite permutations of the phase space (Theorem \ref{main-discr}).
The Radon-Nikodym derivative is computed explicitly and is given by a regularized multiplicative functional.
 Theorem \ref{main-cont} applies, in particular, to
the sine-process,  as well as to determinantal point processes with the Bessel and the Airy kernels; Theorem \ref{main-discr}  to the discrete sine-process and the Gamma kernel process.
The paper  answers a question of Grigori Olshanski.
\end{abstract}

\author{ Alexander I. Bufetov }
\address{Aix-Marseille Universit{\'e}, CNRS, Centrale Marseille, I2M, UMR 7373}
\address{39 rue F. Joliot-Curie, Marseille, France }
\address{Steklov Mathematical Institute of RAS, Moscow, Russia}\address{
The Institute for Information Transmission Problems, Moscow, Russia}\address{
National Research University Higher School of Economics}
\address{Moscow, Russia}
\address{The Chebyshev Laboratory, Saint-Petersburg State University}
\address{ Saint-Petersburg, Russia} 

\email{alexander.bufetov@univ-amu.fr}
\maketitle

\section{Introduction}
\subsection{Outline of the main results.}
G. Olshanski  \cite{GO-Adv} established the quasi-invariance
of the determinantal measure corresponding to the Gamma kernel under the group of finite permutations of ${\mathbb Z}$ and expressed the Radon-Nikodym derivative as a multiplicative functional. S. Ghosh and Y. Peres \cite{GP}, \cite{G} showed, for the Ginibre ensemble
and the Gaussian zero process, that the conditional distribution of particles in a bounded
domain, with the configuration  fixed in the exterior, is equivalent to the Lebesgue measure.

In this paper we take a determinantal point process on ${\mathbb R}$  governed by an orthogonal projection onto a closed subspace $L$
with the following property: given  $p\in {\mathbb R}$ and $\varphi\in L$ satisfying $\varphi(p)=0$, we have
\begin{equation}\label{division-prop}
\frac{\varphi(x)}{x-p}\in L;
\end{equation}
see Assumption \ref{divxp} below for the precise formulation.

In particular, as we check below, Assumption \ref{divxp} holds for projections governed by
kernels
admitting an {\it integrable}  representation
$$
\displaystyle \frac{A(x)B(y) - B(x)A(y)}{x-y}.
$$
The term integrable comes from the connection with the theory of integrable systems discovered by Its, Izergin, Korepin and Slavnov in \cite{IIKS}.

The main results of this paper, Theorems \ref{main-cont} and \ref{main-discr}, establish that, under some additional assumptions,  the measure class of the corresponding
determinantal measures is preserved, in the continuous case,
under the group of diffeomorphisms with compact support (Theorem \ref{main-cont});
in the discrete case, under the group of finite permutations of the phase space (Theorem \ref{main-discr}).
The key step in the proof is the equivalence of reduced Palm measures corresponding to different $l$-tuples of points $(p_1, \dots, p_l)$, $(q_1, \dots, q_l)$   in the phase space; the corresponding Radon-Nikodym derivative  is the regularized  multiplicative functional corresponding to the function
\begin{equation}\label{gpq-cont-1}
\left(\displaystyle \frac{(x-p_1)...(x-p_l)}{(x-q_1)...(x-q_l)}\right)^2.
\end{equation}
The Radon-Nikodym derivative thus has similar form for all the processes with integrable kernels;
the normalizing constants do, of course, depend on the specific process.

 Olshanski \cite{GO-Adv} proves the quasi-invariance of the Gamma-kernel process by a limit transition
from  finite-dimensional approximations. The argument in this paper is direct: first,
it is shown that the Palm subspaces corresponding to conditioning at points $p$ and $q$
 are taken one   to the other by multiplication by the function $(x-p)/(x-q)$; after which, the proof is completed using a general
result of  \cite{Buf-CIRM}, \cite{Buf-umn}  that multiplying  the range of the projection operator $\Pi$ inducing a determinantal measure $\Prob_{\Pi}$
by a function $g$, corresponds, under certain additional assumptions, to taking the product of the
determinantal measure $\Prob_{\Pi}$ by the multiplicative functional $\Psi_g$ induced by the function $g$.
The key technical step is the regularization of divergent multiplicative functionals.

This paper is devoted to determinantal point processes governed by orthogonal projections; in the case of contractions, quasi-invariance is due to Camilier and Decreusefond \cite{cam-decr}; note that in their case the Radon-Nikodym derivative exhibits a much more sensitive dependence on the specific kernel.

\subsection{Projection operators and determinantal point processes}
\subsubsection{Operators and kernels}
Let $\mu$ be a $\sigma$-finite Borel measure on $\mathbb{R}$; for example, $\mu$ can be the Lebesgue measure on $\mathbb R$
or on $\mathbb R_+$ or else the counting measure on ${\mathbb Z}$. The inner product in $L_2({\mathbb R}, \mu)$ will be denoted $\langle, \rangle$.
Let $L \subset L_2(\mathbb{R}, \mu)$ be a closed subspace, and let $\Pi$ be the corresponding operator of orthogonal projection.
We assume that the operator $\Pi$ is locally of trace class and admits a kernel, for which, slightly abusing notation, we keep the same symbol $\Pi$. We let $\Prob_{\Pi}$ be the determinantal measure on the space $\Conf({\mathbb R})$
of configurations on $\mathbb R$ induced by the operator $\Pi$ (see Section 2 below for detailed definitions).
All kernels considered in this paper will always be supposed to satisfy the following
\begin{assumption}\label{ass-ker}
There exists a set $U\subset {\mathbb R}$, satisfying $\mu(\mathbb{R}\setminus U)=0$,
such that
\begin{enumerate}
\item For any $q\in U$ the function
$v_q(x)=\Pi(x,q)$ lies in $L_2(\mathbb{R}, \mu)$ and for any $f\in L_2(\mathbb{R}, \mu)$
we have
$$
\Pi f(q)=\langle f, v_q\rangle.
$$
In particular, all functions in $L$ are defined everywhere on $U$.
\item The diagonal values $\Pi(q,q)$ of the kernel $\Pi$ are defined
for all $q\in U$. We have $\langle v_q, v_q\rangle =\Pi(q,q)$, and, for any
bounded Borel subset $B\subset \mathbb{R}$, we have
$$
\tr (\chi_B\Pi\chi_B)=\displaystyle \int\limits_B \Pi(q,q) d\mu(q).
$$
\item For any $q\in U$ and any $\varphi\in L$ satisfying $\varphi(q)=0$, we have
$$
\frac{\varphi(x)}{x-q}\in L_2(\mathbb{R}, \mu).
$$
\end{enumerate}
\end{assumption}
The first assumptions automatically hold, for instance, for continuous reproducing kernels; the subset $U$ is introduced in order to allow
the consideration of kernels defined on subsets of ${\mathbb R}$ such as, for example, the Bessel kernel defined on ${\mathbb R}_+$
(see Tracy-Widom \cite{TracyWidom}).

The last condition is automatically satisfied once the kernel is sufficiently smooth:
indeed, let $\varphi\in L$ have norm $1$ and be such that $\varphi(q)=0$, let
\begin{equation}\label{ker-orth-q}
{\Pi}^q(x,y)=\Pi(x,y)-\frac{\Pi(x,q)\Pi(q,y)}{\Pi(q,q)}
\end{equation}
be the kernel of the orthogonal projection onto the space $L(q)$, the orthogonal complement of $v_q$ in $L$.
Finally, let $\tilde \Pi$ be the kernel of the orthogonal projection onto the orthogonal complement
of $\varphi$ in $L(q)$.
For any $x\in U$, by definition, we  have
$$
\Pi^qv_x=\langle v_x, \varphi\rangle v_x+\tilde \Pi v_x,
$$
whence, taking the inner  product with $v_x$, we obtain
$$
\Pi^q(x,x)=|\varphi(x)|^2+\langle \tilde \Pi v_x, v_x\rangle.
$$
Using (\ref{ker-orth-q}) and the smoothness of the kernel $\Pi$, we now write the Taylor series for $\Pi^q(x,x)$ in a small neighbourhood of $q$ and obtain $\Pi^q(x,x)=O(|x-q|^2)$, whence also
$|\varphi(x)|=O(|x-q|)$, which implies the desired last condition of Assumption \ref{ass-ker}.

\subsubsection{Palm subspaces}
Given $l$ distinct points $q_1, \ldots, q_l \in \mathbb{R}$, we set
$$
L(q_1, \ldots, q_l) = \{ \varphi \in L : \varphi(q_1) = \ldots = \varphi(q_l) = 0 \};
$$
we denote $\Pi^{q_1, \ldots, q_l}$ the operator of orthogonal projection onto the subspace
$L(q_1, \ldots, q_l)$, and we denote by  $\Prob^{q_1, \dots, q_l}_{\Pi}$ the determinantal measure governed by the projection
 $\Pi^{q_1, \ldots, q_l}$; by the Shirai-Takahashi Theorem,   $\Prob^{q_1, \dots, q_l}_{\Pi}$ is the reduced Palm measure of the
determinantal measure $\Prob_{\Pi}$ with respect to the points $q_1, \dots, q_l$ (the definition and properties of Palm measures are recalled in detail in Section 2 below).

\subsubsection{The main assumption}
We now formulate  our main assumption on our projection operators.

\begin{assumption}\label{divxp}
  If $p\in U$ and $\varphi\in L$ are such that
$\varphi(p)=0$, then there exists $\psi\in L$ such that
\begin{equation}\label{phixpsi}
\varphi=(x-p)\psi.
\end{equation}
\end{assumption}
{\bf Remark.} If the measure $\mu$  does not admit atoms, we can simply say that for $p\in U$ and $\varphi\in L$ satisfying $\varphi(p)=0$  we have (\ref{division-prop}).
In the discrete case, if $\psi$ satisfies (\ref{phixpsi}), then for any $\alpha \in \mathbb {R}$, letting $\delta_p$ be the delta-function at $p$, we see that the function $\psi+\alpha\delta_p$ also satisfies (\ref{phixpsi}). Our requirement is  that one of these functions does belong to $L$.
Induction gives
\begin{proposition}\label{p1-pl}
For any $l\in {\mathbb N}$, if $p_1, \dots, p_l\in U$  are distinct and $\varphi\in L$ satisfies
$\varphi(p_1)=\dots =\varphi(p_l)=0$, then there exists $\psi\in L$ such that
\begin{equation}\label{phixpsi-l}
\varphi=\prod\limits_{i=1}^l\left(x-p_i\right)\psi.
\end{equation}
\end{proposition}

\subsubsection{Change of variables.}

Given a Borel measure $\mu$ on a Borel space $\mathscr X$ and a Borel automorphism $T$ of $\mathscr X$, denote by $\mu\circ T$ the measure defined by $\mu\circ T(\mathscr Z)=
\mu(T(\mathscr Z))$ for all Borel subsets $\mathscr Z\subset \mathscr X$.
Since $T$ is invertible, the measure $\mu\circ T$ is well-defined, and, for any $\mu$-integrable Borel function $f$ on $\mathscr X$, satisfies
$$
\displaystyle \int\limits_{\mathscr X} f\circ T d\mu\circ T= \displaystyle \int\limits_{\mathscr X} f d\mu.
$$
For  a nonnegative Borel function $f$, we also have  $(f\mu)\circ T=(f\circ T)(\mu\circ T)$.

Let $F$ be a Borel automorphism of $\mathbb{R}$ sending bounded sets to bounded sets. The automorphism $F$ acts on $\mathrm{Conf}(\mathbb{R})$ by sending a configuration $X$ to the
 configuration $F(X) = \{ F(x), x \in X \}$;
slightly abusing notation, we keep the same symbol $F$ for this induced action.

Assume additionally that the measure $\mu$ is quasi-invariant under $F$.
Then the measure $\mathbb{P}_\mathcal{\Pi}\circ F^{-1}$ is determinantal with kernel
$$
F^*\Pi(x, y) = \sqrt{\frac{d\mu\circ F}{d\mu}(x) \frac {d\mu\circ F}{d\mu}(y)} \Pi(F(x), F(y)).
$$
The kernel $F^*\Pi$ induces the operator of orthogonal projection onto the subspace
$$
L_{F^*\Pi} = F_* L = \{ \sqrt{F^{\prime}}\cdot \varphi \circ F, \varphi \in L \}.
$$

\subsection{Integrable kernels}
A key particular case is when our kernel $\Pi$ has integrable form
: there exists an open set $U\subset {\mathbb R}$ satisfying $\mu({\mathbb R}\setminus U)=0$ and  linearly independent smooth functions $A$, $B$  defined on $U$ such that
\begin{equation}
\label{piab}
\Pi(x, y) = \displaystyle \frac{A(x)B(y) - A(y)B(x)}{x-y}, x\neq y.
\end{equation}
We  assume that the functions $A,B$ never simultaneously take value $0$ on $U$.
For $p\in U$ we have
$$
v_p(x)=\displaystyle \frac{A(p)B(x) - B(p)A(x)}{p-x};
$$
We have $v_p\in L_2({\mathbb R}, \mu)$ for any $p\in U$ and
for any $\varphi\in L_2({\mathbb R}, \mu)$ we have
$$
\Pi\varphi(p)=\langle \varphi, v_p\rangle.
$$

We consider two cases:
\begin{enumerate}
\item the continuous case: for any $p\in \mathbb R$, $\mu(\{p\})=0$;
\item the discrete case: $\mu$ is the counting measure on a countable subset $E\subset {\mathbb R}$ without accumulation points.
\end{enumerate}
In the continuous case we make the additional requirement
\begin{equation}\label{pi-diag}
\Pi(x, x) = A'(x)B(x)-A(x)B'(x).
\end{equation}
on diagonal values of the kernel $\Pi$;
in the discrete case, when the measure $\mu$ is the counting measure on a countable subset $E\subset {\mathbb R}$ without accumulation points,
the integrability assumption only concerns off-diagonal entries of the kernel $\Pi(x,y)$, and the smoothness assumption is not needed: $A$, $B$ are just arbitrary functions defined on $E$. Note also that the third requirement of Assumption \ref{ass-ker} is only needed in the continuous case.

Note also that the functions $A,B$ in the definition of integrability are not unique: for example, if one makes a linear unimodular change of variable
\begin{equation}\label{chg-var-ab}
(A, B) \to (\alpha_{11}A + \alpha_{12}B, \alpha_{21}A + \alpha_{22}B),  \ \alpha_{11}\alpha_{22} - \alpha_{12}\alpha_{22} = 1,
\end{equation}
then the formula (\ref{piab}) remains valid.

The following proposition plays a crucial r{\^o}le in what follows.
\begin{proposition}\label{intxpred}
An operator of orthogonal projection admitting an integrable kernel satisfies Assumption \ref{divxp}.
\end{proposition}

\subsection{The main result in the continuous case.}

In this subsection we assume that the measure $\mu$ does not admit atoms.
Let $p_1, \ldots, p_l, q_1, \ldots, q_l \in \mathbb{R}$ be distinct.
For $R>0$, $\varepsilon>0$ and a configuration $X$ on $\mathbb{R}$  write
$$
{\overline \Psi}_{R, \varepsilon}(p_1, \ldots, p_l; q_1, \ldots, q_l; X) =
C(R, \varepsilon)\times\prod_{x \in X, |x| \leq R, \min |x-q_i|\geq \varepsilon} \prod_{i=1}^l \left( \displaystyle \frac{x-p_i}{x-q_i} \right)^2,
$$
where the constant $C(R, \varepsilon)$ is chosen in such a way that
\begin{equation}
\displaystyle \int\limits_{\Conf({\mathbb R})}
{\overline \Psi}_{R, \varepsilon}(p_1, \ldots, p_l; q_1, \ldots, q_l; X) d\Prob_{\Pi}^{q_1, \dots, q_l}=1.
\end{equation}
We will  need
 \begin{assumption}\label{xsq}
The kernel $\Pi$ satisfies
\begin{equation}
\displaystyle \int_{\mathbb R} \displaystyle \frac{\Pi(x,x)}{1+x^2}d\mu(x)<+\infty.
\end{equation}
\end{assumption}

\begin{proposition}\label{regmult-cont}
If the kernel $\Pi$ of an orthogonal projection operator satisfies  Assumptions \ref{ass-ker}, \ref{divxp}, \ref{xsq},
then the limit
\begin{equation}\label{defpsireps}
{\overline \Psi}(p_1, \ldots, p_l; q_1, \ldots, q_l; X) = \lim\limits_{R \to \infty, \varepsilon\to 0} \Psi_{R, \varepsilon} (p_1, \ldots, p_l; q_1, \ldots, q_l; X)
\end{equation}
exists  in
$L_1(\mathrm{Conf}(\mathbb{R}), \mathbb{P}_{\Pi^{q_1, \ldots, q_l}})$.
\end{proposition}

{\bf Remark.}  The limit in $R$ and $\varepsilon$ in (\ref{defpsireps}) is understood in the totality of variables, in particular, regardless of order.

\begin{theorem}\label{main-cont}
Let $\mu$ be a continuous measure on ${\mathbb R}$. Let $\Pi$  be a kernel inducing a locally trace-class operator of orthogonal projection in $L_2({\mathbb R}, \mu)$ and satisfying  Assumptions \ref{ass-ker}, \ref{divxp}, \ref{xsq}. Then
\begin{enumerate}
\item For any $l\in {\mathbb N}$ and  two $l$-tuples of distinct points $p_1, \dots, p_l, q_1, \dot, q_l\in U$, we have
$$
\frac {d\Prob_{\Pi}^{p_1, \dots, p_l}}{d\Prob_{\Pi}^{q_1, \dots, q_l}}={\overline \Psi}(p_1, \ldots, p_l; q_1, \ldots, q_l).
$$

\item Let $F: {\mathbb R}\to {\mathbb R}$ be a Borel automorphism acting as the identity beyond a bounded open set $V\subset {\mathbb R}$ and such that $\mu$  is quasi-invariant under $F$. For $\mathbb{P}_{\Pi}$-almost every configuration $X\in \Conf({\mathbb R})$,  $X\bigcap V=\{q_1, \dots, q_l\}$, we have
\begin{multline}
\displaystyle \frac{d\mathbb{P}_{\Pi}\circ F}{d\mathbb{P}_{\Pi}}(X)={\overline \Psi}(F(q_1), \dots, F(q_l); q_1, \dots, q_l; X)\times
 \\ \times \displaystyle \frac{\det({ \Pi}(F(q_i), F(q_j))_{i,j=1,\dots,l}}{\det({\Pi}(q_i,q_j))_{i,j=1,\dots,l}}\times \\ \times \displaystyle \frac{d\mu\circ F}{d\mu}(q_1)\ldots \displaystyle \frac{d\mu\circ F}{d\mu}(q_l).
 \end{multline}
\end{enumerate}
\end{theorem}

{\bf Remark.} The open set $V$ can be chosen in many ways; the resulting value of the Radon-Nikodym derivative is of course the same.

For example, Theorem \ref{main-cont} applies to the sine-process as well as to the Airy and Bessel point processes of Tracy and Widom \cite{TW-Airy}, \cite{TracyWidom}.

{\bf Remark.} In the sequel to this paper, in joint work with Yanqi Qiu we obtain  quasi-invariance results
for determinantal point processes corresponding to Hilbert spaces of holomorphic functions on the plane and on the disc.

\subsection{The main result in the discrete case.}

The main result is similar in the discrete case except that we also need to consider measures  conditional on the absence of particles and that, in order to ensure quasi-invariance of our measures under the infinite symmetric group, we impose the extra restriction that our subspace $L$ not contain functions with finite support.

Let $E\subset {\mathbb R}$ be a countable subset without accumulation points, endowed with the counting measure.
The analogue of Assumption \ref{xsq} in the discrete case is
\begin{assumption} \label{nsq} The subset $E$ satisfies
\begin{equation}
\sum\limits_{n\in E} \displaystyle \frac 1{1+n^2}<+\infty.
\end{equation}
\end{assumption}
Let $\Pi$ be a kernel inducing an operator of orthogonal projection onto a subspace $L\subset L_2(E)$, and let $\Prob_{\Pi}$ be the corresponding determinantal measure on the space $\Conf(E)$ of configurations on $E$. The infinite symmetric group naturally acts on $E$ by finite permutations and induces the corresponding natural action on $\Conf(E)$.
Given $l\in {\mathbb N}$, $m<l$ and an $l$-tuple $(p_1, \dots, p_l)$ of distinct points in $E$ such that there does not exist a nonzero function
in $L$ supported on the set $\{p_1, \dots, p_l\}$,
we introduce a closed subspace  $L(p_1, \dots, p_m, \breve {p}_{m+1}, \dots, \breve {p}_l)$ by the formula
\begin{multline}\label{lppbreve}
  L(p_1, \dots, p_m, \breve {p}_{m+1}, \dots, \breve {p}_l) =\\=
\{  \chi_{E\setminus \{{p}_{m+1}, \dots, {p}_l\}}\varphi : \varphi \in L, \varphi (p_1) =\dots=\varphi(p_m)= 0 \}.
\end{multline}
Let $\Pi^{p_1, \dots, p_m, \breve {p}_{m+1}, \dots, \breve {p}_l}$ be the corresponding orthogonal projection operator,  $\Prob_{\Pi^{p_1, \dots, p_m, \breve {p}_{m+1}, \dots, \breve {p}_l}}$
the corresponding determinantal measure.
The determinantal measure  $\Prob_{\Pi^{p_1, \dots, p_m, \breve {p}_{m+1}, \dots, \breve {p}_l}}$ has the following probabilistic meaning (see Section 2 below for details): consider the conditional measure of $\Prob_{\Pi}$ with respect to the condition that there be particles at
positions $p_1, \dots, p_m$ and holes in positions $p_{m+1}, \dots, p_l$; now {\it remove} the particles at $p_1, \dots, p_m$;
the resulting ``reduced'' conditional measure is precisely $\Prob_{\Pi^{p_1, \dots, p_m, \breve {p}_{m+1}, \dots, \breve {p}_l}}$.

Take $R>0$, $m\leq l$, a permutation $\sigma$ of the points $p_1, \dots, p_l$,
and define
$$
{\overline \Psi}_R(p_1, \dots, p_l, m, \sigma; X)=C_R \prod\limits_{x\in X:|x|\leq R}
\prod\limits_{i=1}^m \left(\displaystyle \frac{x-\sigma(p_i)}{x-p_i}\right)^2\chi_{E\setminus \{p_1, \dots, p_l\}}(x),
$$
where the positive constant $C_R$ is chosen in such a way that
$$
\displaystyle \int\limits_{\Conf(E)} {\overline \Psi}_R(p_1, \dots, p_l, m, \sigma)d\Prob_{\Pi^{p_1, \dots, p_m, {\breve p}_{m+1}, \dots, {\breve p}_l}}=1.
$$

{\bf Remark}. Since the subspace $L$ does not admit functions supported on $\{p_1, \dots, p_l\}$, the set of configurations
having no particles in positions $p_1, \dots, p_l$ has positive probability, consequently, our functional ${\overline \Psi}_R(p_1, \dots, p_l, m, \sigma; X)$ is positive with positive probability.

\begin{proposition}\label{regmult-discr}
Let $E$ be a countable subset of $\mathbb R$ without accumulation points
satisfying  Assumption \ref{nsq}.
Let $\Pi$ be a kernel inducing an operator of orthogonal projection on $L_2(E)$. Let $p_1, \dots, p_l\in E$ be distinct points such that there does not exist a nonzero function
in $L$ supported on the set $\{p_1, \dots, p_l\}$.
Then, for any  $m\leq l$ and any permutation $\sigma$ of $p_1, \dots, p_l$, the limit
$$
{\overline \Psi}(p_1, \dots, p_l, m, \sigma)=\lim\limits_{R\to\infty}
{\overline \Psi}_R(p_1, \dots, p_l, m, \sigma)
$$
exists in $L_1(\Conf(E), \Prob_{\Pi^{p_1, \dots, p_m, \breve p_{m+1}, \dots, \breve p_l}})$.
\end{proposition}

Let ${\bf C}(p_1, \dots, p_m, \breve{p}_{m+1}, \dots, \breve{p}_l)$ be the set of configurations on $E$ containing
exactly one particle in each of the positions $p_1, \dots, p_m$ and no particles in the positions $p_{m+1}, \dots, p_l$.

We are now ready to formulate the main result in the discrete case,  the quasi-invariance of determinantal measures with integrable kernels  under the natural action of the infinite symmetric group on $\Conf(E)$.
Given a permutation $\sigma$ of points $p_1, \dots, p_l$ of the set $E$, slightly abusing notation, we use the same symbol both for the bijection of $E$ that acts as $\sigma$ on $\{p_1, \dots, p_l\}$ and as the identity on $E\setminus \{p_1, \dots, p_l\}$ and the automorphism induced by this bijection on the space $\Conf(E)$ of configurations on $E$.
\begin{theorem}\label{main-discr} Let  $E$ be a countable subset of $\mathbb R$ without accumulation points satisfying (\ref{nsq}).
Let $\Pi$ be a kernel inducing an operator of orthogonal projection onto a closed subspace  $L\subset L_2(E)$ and satisfying  Assumptions \ref{ass-ker}, \ref{divxp}.
Let $p_1, \dots, p_l$ be distinct points in $E$ such that there does not exist a nonzero function
in $L$ supported on the set $\{p_1, \dots, p_l\}$.
 Then for any $m\leq l$, any permutation $\sigma$ of the points $p_1, \dots, p_l$ and $\Prob_{\Pi}$-almost every $X\in {\bf C}(p_1, \dots, p_m, \breve{p}_{m+1}, \dots, \breve{p}_l)$, we have
\begin{equation}
\displaystyle \frac{d\mathbb{P}_{\Pi}\circ \sigma}{d\mathbb{P}_{\Pi}}(X)={\overline \Psi}(p_1, \dots, p_l, m,  \sigma; X) \times \displaystyle \frac{\det\left({ \Pi}(\sigma(p_i), \sigma(p_j))\right)_{i,j=1,\dots,m}}{\det\left({\Pi}(p_i,p_j)\right)_{i,j=1,\dots,m}}.
 \end{equation}
In particular, if the subspace $L$ does not contain functions supported on finite sets, then the measure
$\Prob_{\Pi}$ is quasi-invariant under the natural action of the infinite symmetric group on $\Conf(E)$.
\end{theorem}

For example, the discrete sine-process of Borodin, Okounkov and Olshanski \cite{BOO} as well as the Gamma kernel process of Borodin and Olshanski \cite{BO-gamma} satisfy all the assumptions of Theorem \ref{main-discr}.

{\bf Remark.} By the Theorem of Ghosh \cite{Ghosh}, the sine-process, discrete or continuous,  is {\it rigid}: if, for  a bounded subset $B $ and a configuration $X$, we let $ \#_B(X)$ stand for the number of particles of $X$ lying in $B$, and, for any Borel subset $C$ we let $\mathcal{F}_C$ be the $\sigma$-algebra generated by all random variables of the form $\#_B, B\subset C$, then, for any bounded $B$, the random variable $\#_B$ is measurable with respect to the completion, under the sine-process,
of the sigma-algebra $\mathcal{ F}_{B^c}$, where $B^c$ stands for the complement of $B$.
 As R. ~Lyons, developing the method of \cite{BLPS}, showed in Theorem 7.15  of \cite{Lyons}, the tail sigma-algebra of the discrete sine-process is trivial. It follows now that the {\it symmetric} sigma-algebra of the sine-process is trivial as well: in other words,
the discrete sine-process is ergodic with respect to the action of the infinite symmetric group.
This argument holds, of course, for any rigid point process.

To further illustrate Theorem \ref{main-discr}, we now write the Radon-Nikodym derivative for a transposition of two points $p,q\in E$.
	Set
$$
	L(p,\breve{q}) = \{ \chi_{E\backslash \{p,q\}} \varphi, \text{ } \varphi \in L,  \varphi(p)=0\}.
$$
	and let $\mathbb{P}_{\Pi}^{p, \breve{q}}$ be the determinantal measure corresponding to the operator of orthogonal projection onto the subspace $L (p, \breve{q})$.
The subspace $L (q, \breve{p})$ and the measure $\mathbb{P}_{\Pi}^{q, \breve{p}}$ are defined in the same way.
	Write
$$
	\overline{\Psi}_N (p, q; X) = C_{p, q} \times \prod\limits_{ x \in X, |x| \le N } {\left( \displaystyle \frac{x-p}{x-q} \right)}^2,
$$
 	where the constant $C_{p,q}$ is chosen in such a way that
$$
	\displaystyle \int\limits_{\Conf(E)} \overline{\Psi}_N (p, q; X) \text{ }d\mathbb{P}_{\Pi}^{p, \breve{q}}(X) = 1.
$$
	By definition, $\mathbb{P}_{\Pi}^{p, \breve{q}}$-almost all configurations $X$ on $E$ contain no particles either at $p$ or at $q$,  so the function $\overline{\Psi}_N$ is well-defined; by definition it is bounded.

\begin{proposition}
	The limit $\overline{\Psi}(p, q; X) = \lim\limits_{N \to \infty} \overline{\Psi}_N (p, q; X)$ exists  in $L_1(\Conf(E), \mathbb{P}_{\Pi}^{p, \breve{q}})$.
\end{proposition}
	The Radon-Nikodym derivative of $\Prob_{\Pi}$ under the action of the permutation $\sigma_{pq}$ is now given by the following
\begin{proposition}\label{discrete-sine}	For $\mathbb{P}_{\Pi}$-almost all $X\in \Conf(E)$ the following holds.

If $p\notin X, q\in X$, then
$$		
	\displaystyle \frac{d\mathbb{P}_{\Pi} \circ \sigma_{pq}}{d\mathbb{P}_{\Pi}}(X) =
\overline{\Psi}(p, q; X)\cdot \frac{\Pi(p,p)}{\Pi(q,q)}.
$$

If $p\in X, q\notin X$, then
$$
	\displaystyle \frac{d\mathbb{P}_{\Pi} \circ \sigma_{pq}}{d\mathbb{P}_{\Pi}}(X) = \overline{\Psi}(q, p; X)
\cdot \frac{\Pi(q,q)}{\Pi(p,p)}.
$$

If $p,q\in X$ or $p,q\notin X$, then
$$
	\displaystyle \frac{d\mathbb{P}_{\Pi} \circ \sigma_{pq}}{d\mathbb{P}_{\Pi}}(X) = 1.
$$
\end{proposition}

{\bf Remark.} If $E$ is a countable set,  $\Prob$  a Gibbs  measure on $\Conf(E)$ corresponding to
the potential $U$ of pairwise interaction of particles (cf. e.g. Sinai \cite{sinai}), $p,q$ are points in $E$ and $\sigma_{pq}$ the transposition of $p$ and $q$, then, for almost every configuration $X$, conditioned to contain a particle at $q$ but not at $p$,  by definition, we have
$$
\frac{d\Prob\circ \sigma_{pq} }{d\Prob}(X)=\prod\limits_{x\in X: x\neq q} \exp(U(p,x)-U(q,x)).
$$
The quasi-invariance property established in this paper is  the analogue of the Gibbs property  for determinantal point processes governed by integrable kernels.  The potential is $U(x,y)=2\log|x-y|$.

\subsection{Examples of regularized multiplicative functionals.}

Regularization of a multiplicative functional can take different form depending on the specific process.
We illustrate this by two examples.

{\it The Sine-Process.} The argument below is valid for the continuous sine-process as well as the discrete sine-process. The sine-process is stationary, therefore, for almost every configuration $X$ the series
\begin{equation}\label{sum1x}
\sum\limits_{x\in X: x\neq 0} \frac 1x
\end{equation}
{\it diverges} absolutely since so does the harmonic series. Nonetheless, the series (\ref{sum1x}) converges {\it conditionally} in principal value: the limit
$$
\lim\limits_{N\to\infty} \sum\limits_{x\in X: x\neq 0, |x|\leq N} \frac 1x
$$
is almost surely finite and, as we shall check below, has finite variance.
Similarly,  for distinct points $p_1, \ldots, p_l, q_1, \ldots, q_l$, taken in ${\mathbb R}$ in the continuous case and in ${\mathbb Z}$ in the discrete case, the limit
\begin{equation}\label{lim-prod-sine}
\lim\limits_{N\to\infty}\prod_{x \in X, |x| \leq N, x\neq q_1, \dots, q_l} \prod_{i=1}^l \left( \displaystyle \frac{x-p_i}{x-q_i} \right)^2,
\end{equation}
 exists and has finite expectation. The normalized mutliplicative functional is in this case precisely
the limit (\ref{lim-prod-sine}) normalized to have expectation $1$.

{\it The Determinantal Point Process with the Gamma-Kernel.}
The determinantal point process with the Gamma-kernel, introduced by Borodin and Olshanski
in \cite{BO-gamma} and for which the quasi-invariance under the action of the infinite symmetric group is due to Olshanski \cite{GO-Adv}, is a point process on the phase space ${\mathbb Z}^{\prime}=1/2+{\mathbb Z}$ of half-integers  such that for almost every configuration $X$ we have
\begin{equation}\label{sum-part-hole}
\sum\limits_{x\in X: x>0} \frac 1x<+\infty, \ \sum\limits_{y\notin X: y<0} \frac 1{|y|}<+\infty.
\end{equation}
Furthermore, each sum in (\ref{sum-part-hole}), considered as a random variable on the space of configurations on ${\mathbb Z}^{\prime}$, has finite variance with respect to the determinantal point process with the Gamma-kernel.

For $p,q\in {\mathbb Z}^{\prime}$, the normalized multiplicative functional corresponding to the function
$g(x)=((x-p)/(x-q))^2$ will therefore have the form
$$
C\cdot \prod\limits _{x\in X, x>0} g(x)  \cdot \prod\limits _{y\notin X: y<0} g^{-1}(y),
$$
where the constant $C$ is chosen in such a way that the resulting expression have expectation $1$.

\subsection{Outline of the argument}
	
	We start with the discrete case and illustrate the argument in the specific case of a transposition of two distinct points $p,q\in E$.  A theorem due to Lyons~\cite{Lyons}, Shirai-Takahashi~\cite{ShirTaka2} states that the measure $\mathbb{P}_{\Pi}^{p, \breve{q}}$ is the conditional measure of $\mathbb{P}_{\Pi}$ on the subset of configurations containing a particle at $p$ and not containing a particle at $q$.

	{\it Step 1. The Relation Between Palm Subspaces.} The key point in the proof of Proposition \ref{discrete-sine} is the equality
\begin{equation}\label{subsp-sine1}
	L(p, \breve{q}) = \displaystyle \frac{x-p}{x-q} L(q, \breve{p}),
\end{equation}
which it is more convenient to rewrite in the form
\begin{equation}\label{subsp-sine2}
	L(p, \breve{q}) = \chi_{E\setminus\{p, q\}} \displaystyle \frac{x-p}{x-q} L(q, \breve{p}).
\end{equation}

	The equality (\ref{subsp-sine2}) directly follows from the {\it integrability} of the discrete sine-kernel. The remainder of the argument shows that the relation (\ref{subsp-sine2}) implies the relation
\begin{equation}\label{pq-eq}
	\mathbb{P}_{\Pi}^{p, \breve{q}} = \overline{\Psi}(p, q) \mathbb{P}_{\Pi}^{q, \breve{p}},
\end{equation}
	which, in turn, is a reformulation of Proposition \ref{discrete-sine}.

	{\it Step 2. Multiplicative functionals of determinantal point processes.}	
	Given a function $g$ on $\mathbb{Z}$, the multiplicative functional $\Psi_g$ is defined on $\Conf(E)$ by the formula
$$	
	\Psi_g (X) = \prod\limits_{x \in X} g(x).
$$
provided that the infinite product in the right-hand side converges absolutely.
	
	 At the centre of the argument lies the result of ~\cite{Buf-umn} that can informally be summarized as follows: a determinantal measure times a multiplicative functional is, after normalization, again a determinantal measure. More precisely, let $g$ be a positive function on $E$ bounded away from 0 and $\infty$, and let $\Pi$ be an operator of orthogonal projection in $L_2(E)$ onto a closed subspace $L$. Let $\Pi^g$ be the operator of orthogonal projection onto the subspace $\sqrt{g}L$. Then, under certain additional assumptions we have
\begin{equation}\label{mult-1}
	\mathbb{P}_{\Pi^g} = \displaystyle \frac{\Psi_g \mathbb{P}_\Pi}{\displaystyle \int\limits_{\Conf(E)} \Psi_g \text{ }d\mathbb{P}_\Pi}
\end{equation}				
 	The relation (\ref{mult-1}) together with the relation (\ref{subsp-sine2}) suggests that the measures $\mathbb{P}_{\Pi}^{p, \breve{q}}$ and $\mathbb{P}_{\Pi}^{q, \breve{p}}$ are equivalent, and the Radon-Nikodym derivative is given by the normalized multiplicative functional corresponding to the function
$$
	g^{p,\breve{q}}(x) = \displaystyle \frac{x-p}{x-q} \chi_{E\setminus \{p,q\}}.
$$

{\it Step 3. Regularization of multiplicative functionals.}
	A technical difficulty arises that in many examples the multiplicative functional corresponding to the function $g^{p,\breve{q}}$   fails to converge absolutely with respect to the measure $\mathbb{P}_{\Pi}^{p, \breve{q}}$; indeed, in many examples (in particular, for stationary determinantal processes on ${\mathbb Z}$), we  have
$$
	\sum\limits_{x\in E} \left|  g^{p,\breve{q}}(x) - 1 \right|  = +\infty
$$
and, consequently,  also
$$
	\sum\limits_{x \in E} \left|  g^{p,\breve{q}}(x) - 1 \right| \cdot {\Pi}^{p,\breve{q}}(x,x) = +\infty.
$$

	In order to resolve this difficulty, we go back to the formula (\ref{mult-1}). For multiplicative functional $\Psi_g$ integrable with respect to a determinantal measure $\mathbb{P}_\Pi$ set
\begin{equation}\label{reg-mult-def}
	\overline{\Psi}_g = \displaystyle \frac{\Psi_g}{\displaystyle \int \Psi_g d\mathbb{P}_\Pi}.
\end{equation}
	The functional $\overline{\Psi}_g$ will be called the normalized multiplicative functional corresponding to $\Psi_g$ and $\mathbb{P}_{\Pi}$. To keep notation lighter, we do not explicitly indicate dependence on $\Pi$; in what follows, the precise measure, with respect to which normalization is taken, will be clear from the context. We now rewrite (\ref{mult-1}) in the form

\begin{equation}\label{mult-2}
	\mathbb{P}_{\Pi^g} = \overline{\Psi}_g \cdot \mathbb{P}_{\Pi}.
\end{equation}
 	The key observation for the remainder of the argument is that the definition of the normalized multiplicative functional $\overline{\Psi}_g$ can be extended in such a way that (\ref{mult-2}) continues to hold for a wider class of functions $g$, for which the multiplicative functional itself diverges almost surely.

	We first explain the idea of this extension for additive functionals. Given a measurable function $f$ on $E$, the corresponding additive functional on $\Conf(E)$ is defined by the formula
$$
	S_f(X) = \sum\limits_{x \in X} f(x)
$$
provided the series in the right hand side converges absolutely. The expectation of the additive functional with respect to $\mathbb{P}_{\Pi}$ is given by the formula

\begin{equation}\label{sin-exp}
	\mathbb{E}_{\mathbb{P}_{\Pi}} S_f = \displaystyle \sum\limits_{x\in E} f(x)\Pi(x,x),
\end{equation}
provided, again, that the series in the right hand side converges absolutely. For the variance of the additive functional we have

\begin{equation}\label{sin-var}
	\mathrm{Var}_{\mathbb{P}_{\Pi}} S_f = \frac12\sum\limits_{x,y \in \mathbb{Z}} (f(x) - f(y))^2 ({\Pi}(x,y))^2.
\end{equation}

	Let $$\overline{S}_f = S_f - \mathbb{E}_{\mathbb{P}_{\Pi}} S_f $$ be the normalized additive functional corresponding to the function $f$.
	It is easy to give examples of functions $f$ for which the sum in the right hand side of (\ref{sin-exp}) diverges while the sum in the right hand side of (\ref{sin-var}) converges.
	For such functions, convergence of the sum in the right hand side of (\ref{sin-var}) allows one to define the normalized additive functional $\overline{S}_f$ by continuity, even though the additive functional $S_f$ itself is not defined.
	In a similar way, for a function $g$ bounded away from 0 and $\infty$ and satisfying
$$
\sum_{x \in E} \left| g(x) - 1 \right|^2 \Pi(x,x) < +\infty,
$$
one can define the normalized multiplicative functional $\overline{\Psi}_g$ even when the multiplicative functional $\Psi_g$ itself is not  defined; the relation (\ref{mult-2}) still holds.
	
 We next check that the normalized multiplicative functional $\overline{\Psi}_{g^{p,q}}$ can be defined with respect to the measure $\mathbb{P}_{\Pi}^{q, \breve{p}}$; note that the function $g^{p,q}$  assumes value zero at $p$, and here we must use the assumption that our subspace $L$ does not contain functions with finite support.
The relation (\ref{subsp-sine2}) now implies the  equality (\ref{pq-eq}), and the outline of the proof of Theorem \ref{main-discr} is complete.

	The proof in the continuous case follows a similar scheme.
	The r{\^o}le of conditional measures is played by reduced Palm measures.  The reduced Palm measure
$\mathbb{P}_{\Pi}^{q_1, \dots, q_l} $
of the measure $\mathbb{P}_{\Pi}$ with respect to $l$ distinct points $q_1,..., q_l \in \mathbb{R}$ is the determinantal measure corresponding to the operator ${\Pi}^{q_1,...,q_l}$ of the orthogonal projection onto the subspace
$$
	L(q_1,...,q_l) = \{ \varphi \in L : \varphi(q_1) = ... = \varphi(q_l) = 0 \}.
$$
	 The continuous analogue of the equality (\ref{subsp-sine2}) is the relation
\begin{equation}\label{subsp-sine}
	L(p_1,...,p_l) = \displaystyle \frac{(x-p_1)...(x-p_l)}{(x-q_1)...(x-q_l)}L(q_1,...,q_l)
\end{equation}
	 valid for $\mu^{\otimes l}$-almost any two $l$-tuples of distinct points $(p_1,...,p_l)$, $(q_1,...,q_l)$ in $\mathbb R$.

	The next step is to regularize the multiplicative functional corresponding to the function
\begin{equation}\label{gpq-cont}
\displaystyle\displaystyle \frac{(x-p_1)...(x-p_l)}{(x-q_1)...(x-q_l)};
\end{equation}
while the overall scheme of regularization is the same  as in the discrete case,
 additional estimates are needed here because the function (\ref{gpq-cont}) is bounded away neither from zero nor from infinity.

The resulting normalized multiplicative functional $\overline{\Psi}(p_1,...,p_l,q_1,...,q_l)$ is then seen to be the Radon-Nikodym derivative of the reduced Palm measures
$\mathbb{P}_{\Pi}^{p_1, \dots, p_l} $ and $\mathbb{P}_{\Pi}^{q_1, \dots, q_l}$, which, in turn, implies
Theorem \ref{main-cont}.

\subsection{Organization of the paper}
The paper is organized as follows. In Section 2, we collect necessary facts about determinantal point processes, their multiplicative functionals and their Palm measures. We recall the results of \cite{Buf-umn} (see also \cite{Buf-CIRM}) showing that the product  of a determinantal measure with a multiplicative functional is, after normalization,  again a determinantal measure, whose kernel is found explicitly. We  also check that equivalence of reduced Palm measures corresponding to distinct $l$-tuples  of points implies the quasi-invariance of the point process under Borel automorphisms preserving the class of its correlation measures and acting by the identity beyond a bounded set. In Section 3, we start by showing that reduced Palm measures of determinantal point processes given by projection operators with integrable kernels are themselves determinantal point processes given by projection operators with integrable kernels and proceed to verify the key relations  (\ref{palm-mult-cont})
and (\ref{palmeqdisc}) showing that the ranges of projection operators corresponding to reduced Palm measures at distinct points differ by multiplication by a function.

 Proposition \ref{multreg-main} in Section 4 describes the  regularization of  multiplicative functionals. Relations  (\ref{palm-mult-cont})
and (\ref{palmeqdisc}) are then seen to imply that
the reduced Palm measures themselves are equivalent, and that the corresponding  Radon-Nikodym derivative is a regularized multiplicative functional, thus completing the proof of the main Theorems \ref{main-cont} and \ref{main-discr}.

{\bf {Acknowledgements.}}
Grigori Olshanski posed the problem to me and suggested that the  Radon-Nikodym derivative be given
by a multiplicative functional; I am greatly indebted to him. I am deeply grateful to Alexei Klimenko and Cosme Louart for useful discussions. I am deeply grateful to the anonymous referees whose suggestions have
been of great help in improving the presentation.

This work has been carried out thanks to the support of the A*MIDEX project (no. ANR-11-IDEX-0001-02) funded by the  programme ``Investissements d'Avenir '' of the Government of the French Republic, managed by the French National Research Agency (ANR).
This project has received funding from the European Research Council (ERC) under the European Union's Horizon 2020 research and innovation programme (grant agreement No 647133 (ICHAOS)) and has also been funded by Grant MD 5991.2016.1 of the President of the Russian Federation,
by the Russian Academic Excellence Project `5-100' and by the Chaire Gabriel Lam{\'e} at the Chebyshev Laboratory.

Part of this work was done while I was visiting  Institut Henri Poincar{\'e} in Paris and the
Max Planck Institute in Bonn; I am deeply grateful to these institutions for their warm hospitality.

{\bf Remark.} After this work was completed, I became aware of the preprint {\it Absolute continuity and singularity of Palm measures of the Ginibre point process}, arXiv:1406.3913, by Hirofumi Osada and Tomoyuki Shirai, in which, for the Ginibre ensemble, using its finite-dimensional approximations by orthogonal polynomial ensembles, the authors establish the equivalence of reduced Palm measures, conditioned at distinct $l$-tuples of points in ${\mathbb C}$, and represent the Radon-Nikodym derivative as a regularized multiplicative functional.

\section{Point Processes and Palm Distributions.}
\subsection{Spaces of configurations.}
Let $E$ be a locally compact complete metric space.
A {\it configuration} on $E$ is a collection  of points in $E$, called {\it particles},  considered without regard to order and
subject to  the additional requirement that every bounded set contain only finitely many particles of a configuration.
Let $\Conf(E)$ be the space of configurations on $E$.
To a configuration $X\in \mathrm{Conf}(E)$ assign a Radon measure
$
 \sum\limits_{x\in X}\delta_x
$
on the space $E$; this correspondence identifies the space
 $\mathrm{Conf}(E)$ with  the space of integer-valued
Radon measures on $E$.  The
space $\mathrm{Conf}(E)$ is thus endowed with a natural structure of a complete
separable metric space. The Borel structure on the space $\Conf(E)$ can equivalently be defined without
introducing a topology explicitly: namely, for a bounded Borel set $B\subset E$,   let
$$\#_B\colon\Conf(E)\to\mathbb{N}\cup\{0\}$$ be
the function that to a configuration
assigns the number of its particles
belonging to~$B$. The random variables $\#_B$ over all  bounded Borel sets $B\subset E$
determine the  Borel sigma-algebra on $\Conf(E)$.

\subsection{Multiplicative functionals}
We next recall the definition of {\it multiplicative functionals} on spaces of configurations.
Let $g$ be a non-negative measurable function on $E$, and introduce the
{\it multiplicative functional} $\Psi_g:\Conf(E)\to\mathbb{R}$ by the formula
 \begin{equation} \label{mult-fun-def}
 \Psi_g(X)=\prod\limits_{x\in X}g(x).
 \end{equation}
If the infinite product
$\prod\limits_{x\in X}g(x)$ absolutely converges to $0$ or to $\infty$, then we set, respectively,
$\Psi_g(X)=0$ or $\Psi_g(X)=\infty$. If the product in the right-hand side fails to converge absolutely,
then the multiplicative functional is not defined.

\subsection{Point processes.}
A Borel probability measure $\Prob$ on $\Conf(E)$ is called {\it a point process} with phase space $E$.

We recall that the process $\Prob$ is said to admit correlation functions
of order $l$ if for any continuous compactly supported function $f$ on $E^l$
the functional
$$
\sum\limits_{x_1, \dots, x_l\in X} f(x_1, \dots, x_l)
$$
is $\Prob$-integrable; here the sum is taken over all $l$-tuples of distinct particles in $X$. The $l$-th correlation measure $\rho_l$ of the point process $\Prob$ is then
defined by the formula
$$
\ee_{\Prob} \left(\sum\limits_{x_1, \dots, x_l\in X} f(x_1, \dots, x_l)\right)=
\displaystyle \int\limits_{E^l} f(q_1, \dots, q_l)d\rho_l(q_1, \dots, q_l).
$$

By definition, a  point process $\Prob$ is uniquely determined by prescribing joint distributions, with respect to $\Prob$,  of
random variables $\#_{B_1}, \dots, \#_{B_l}$ over all finite collections of disjoint bounded Borel subsets $B_1, \dots, B_l\subset E$.
Since, for arbitrary nonzero complex numbers $z_1, \dots, z_l$ inside the unit circle, the function
$
\prod\limits_{k=1}^l z_k^{\#_{B_k}}
$
is a well-defined multiplicative functional on $\Conf(E)$,  that, moreover, takes values inside the unit circle, a point process $\Prob$ on $\Conf(E)$ is also uniquely determined by prescribing the values of expectations of  multiplicative functionals of this form.

\subsection{Campbell Measures.}

Following Kallenberg \cite{Kallenberg} and Daley--Vere-Jones \cite{DVJ}, we now recall the definition of Campbell measures of point processes.

Take a Borel probability measure $\Prob$ on $\Conf(E)$ of {\it finite local intensity}, that is, admitting the first correlation measure $\rho_1$, or, equivalently, such that
 for any bounded Borel set $B$, the function $\#_B$ is integrable with respect to $\Prob$.
For any bounded Borel set $B\subset E$, by definition we then have
$$
\rho_1(B)=\displaystyle \int_{\Conf(E)} \#_B(X)d\Prob(X).
$$

The {\it Campbell measure}  ${\EuScript C}_{\Prob}$ of a  Borel probability measure $\Prob$ of finite local intensity on $\Conf(E)$
is a sigma-finite measure on $E\times \Conf(E)$ such that for any Borel subsets
$B\subset E$, ${\mathscr Z}\subset \Conf(E)$ we have
$$
{\EuScript C}_{\Prob}(B\times {\mathscr Z})=\displaystyle \int\limits_{{\mathscr Z}} \#_B(X)d\Prob(X).
$$

For a point process admitting correlation functions of order $l$ one can also define the $l$-th iterated Campbell measure ${\EuScript{C}}^{(l)}$  of the point process $\mathbb{P}$, that is, by definition, a measure on $E^l \times \mathrm{Conf}(E)$ such that for any disjoint bounded sets $B_1,...,B_l \subset E$ and any measurable subset $\mathscr Z \subset \mathrm{Conf}(E)$ we have
\begin{equation}\label{camp-def}
{\EuScript C}^{(l)}(B_1 \times ... \times B_l \times \mathscr Z) = \displaystyle \int\limits_{\mathscr Z}
{\#}_{B_1}(X)\times...\times {\#}_{B_l}(X) d\mathbb{P}(X).
\end{equation}

\subsection{Palm Distributions.}
Following Kallenberg \cite{Kallenberg} and Daley--Vere-Jones \cite{DVJ}, we now recall the construction of Palm distributions from Campbell measures.
For a fixed Borel ${\mathscr Z}\subset \Conf(E)$ the Campbell measure ${\EuScript C}_{\Prob}$ induces a sigma-finite measure ${\EuScript C}_{\Prob}^{{\mathscr Z}}$
on $E$ by the formula
$$
{\EuScript C}_{\Prob}^{{\mathscr Z}}(B)={\EuScript C}_{\Prob}(B\times {\mathscr Z}).
$$

By definition, for any Borel subset ${\mathscr Z}\subset \Conf(E)$ the measure
${\EuScript C}_{\Prob}^{{\mathscr Z}}$ is absolutely continuous
with respect to $\rho_1$.
We now take $q\in E$, vary ${\mathscr Z}$ and set
$$
{\hat \Prob}^q({\mathscr Z})=\displaystyle \frac{d{\EuScript C}_{\Prob}^{{\mathscr Z}}}{d\rho_1}\left(q\right).
$$
For  $\rho_1$-almost every $q\in E$, the expression ${\hat \Prob}^q({\mathscr Z})$, considered as a function of ${\mathscr Z}$ (while $q$ stays fixed),
defines a probability measure ${\hat \Prob}^q$ on $\Conf(E)$, {\it the Palm measure} of $\Prob$ at the point $q$.
Equivalently, the Palm measure ${\hat \Prob}^q$ is the canonical conditional measure, in the sense of Rohlin \cite{Rohmes}, of the Campbell measure ${\mathcal C}_{\Prob}$
with respect to the measurable partition of the space $E\times \Conf(E)$
into subsets of the form $\{q\}\times \Conf(E)$, $q\in E$.

Similarly, using iterated Campbell measures one defines iterated Palm measures:
for a fixed Borel ${\mathscr Z}\subset \Conf(E)$ the $l$-th iterated Campbell measure ${\EuScript C}^l_{\Prob}$ induces a sigma-finite measure ${\EuScript C}_{\Prob}^{l,{\mathscr Z}}$
on $E$ by the formula
$$
{\EuScript C}_{\Prob}^{l,{\mathscr Z}}(B)={\EuScript C}_{\Prob}(B\times {\mathscr Z}).
$$
By definition, for any Borel subset ${\mathscr Z}\subset \Conf(E)$ the measure
${\EuScript C}_{\Prob}^{l, {\mathscr Z}}$ is absolutely continuous
with respect to the $l$-th correlation measure $\rho_l$ of our point  process $\Prob$. For  $\rho_l$-almost all $(q_1, \dots, q_l)\in E^l$,
one can therefore define a probability measure ${\hat \Prob}^{q_1, \dots, q_l}$ on $\Conf(E)$
by the formula
$$
{\hat \Prob}^{q_1, \dots, q_l}({\mathscr Z})=\displaystyle \frac{d{\EuScript C}_{\Prob}^{{\mathscr Z}}}{d\rho_1}\left(q_1, \dots, q_l\right).
$$
The measure ${\hat \Prob}^{q_1, \dots, q_l}$ is called
the $l$-th iterated Palm measure of the point process $\Prob$.
The iterated Palm measure ${\hat \Prob}^q$ is the canonical conditional measure, in the sense of Rohlin \cite{Rohmes}, of the Campbell measure ${\mathcal C}^l_{\Prob}$
with respect to the measurable partition of the space $E^l\times \Conf(E)$
into subsets of the form $\{q_1, \dots, q_l\}\times \Conf(E)$, with
$q_1, \dots, q_l\in E$ distinct.

For distinct  points $q_1, \dots, q_l$, the $l$-th iterated
Palm measure of course satisfies
$$
{\hat \Prob}^{q_1, \dots, q_l}=\left(\dots \left({\hat \Prob}^{q_1}\right)^{q_2}\dots\right)^{q_l}.
$$
\subsection{Reduced Palm measures}

By definition, the Palm measure ${\hat \Prob}^{q_1, \dots, q_l}$ is supported on the subset of configurations containing a particle at each position $q_1, \dots, q_l$.
It is often convenient to remove these particles and to define the {\it reduced}
Palm measure $\Prob^{q_1, \dots, q_l}$ as the push-forward of the Palm measure
${\hat \Prob}^{q_1, \dots, q_l}$ under the erasing map
$X\to X\setminus \{q_1, \dots, q_l\}$.
Reduced Palm measures allow one to give a convenient representation for measures of cylinder sets.
Take $X_0\in\Conf(E)$ and $q_1^{(0)}, \dots, q_l^{(0)}\in X_0$.

Take disjoint bounded open sets $B^{(1)},...,B^{(l)}\subset E$ such that $q_i^{(0)}\in  { B^{(i)}}$,
$q_i^{(0)}=X_0\cap {\overline B^{(i)}}$ for all $i=1, \dots, l$.
Set $B=\cup B^{(i)}$ and
take an open set $V\subset E$ with bounded complement,
disjoint from all $B^{(i)}$, and  satisfying  $X_0\setminus \{q_1^{(0)}, \dots, q_l^{(0)}\}\subset V$. Let $\mathscr W$ be a neighbourhood of $X_0\setminus \{q_1^{(0)}, \dots, q_l^{(0)}\}$ in $\Conf(E)$
satisfying $$\mathscr W\subset \{X\in\Conf(E): X\subset V\}.$$
Introduce a neighbourhood $\mathscr Z$ of $X_0$ by setting
\begin{equation}\label{bnsp}
\mathscr Z = \{ X \in\Conf(E):  {\#}_{B^{(1)}}(X)=\dots={\#}_{B^{(l)}}(X)=1, X|_{E\setminus B}\subset  {\mathscr W}\}.
\end{equation}
\begin{proposition}\label{palm-cyl}
We have
$$
	\mathbb{P}(\mathscr Z) = \displaystyle \int\limits_{B^{(1)} \times ... \times B^{(l)}} {{\mathbb{P}}}^{q_1,...,q_l}(\mathscr W)
	d\rho_l(q_1,...,q_l).
$$
\end{proposition}
\begin{proof}
We disintegrate  ${\EuScript{C}}^{(l)}(B^{(1)} \times ... \times B^{(l)} \times \mathscr Z)$ in two  ways.
By definition of iterated Palm measures, we have
$$
	{\EuScript{C}}^{(l)}(B^{(1)} \times ... \times B^{(l)} \times \mathscr Z) = \displaystyle \int\limits_{B^{(1)} \times ... \times B^{(l)}}
	{\hat {\mathbb{P}}}^{q_1,...,q_l}(\mathscr Z) d\rho_l(q_1,...,q_l).$$
By definition of iterated Campbell measures, see (\ref{camp-def}), (\ref{bnsp}), we have
$$
	{\EuScript{C}}^{(l)}(B^{(1)} \times ... \times B^{(l)} \times \mathscr Z) =
\displaystyle \int\limits_{\mathscr Z}
{\#}_{B^{(1)}}(X)\times...\times {\#}_{B^{(l)}}(X) d\mathbb{P}(X)= \mathbb{P}(\mathscr Z).
$$
Passing to reduced Palm measures gives
$
{\hat {\mathbb{P}}}^{q_1,...,q_l}(\mathscr Z)= 	{\mathbb{P}}^{q_1,...,q_l}(\mathscr W),
$
and (\ref{bnsp}) follows. 
\end{proof}

\subsection{Locally  trace class operators and their kernels.}
Let $\mu$ be a sigma-finite Borel measure on $E$. The inner product in $L_2(E, \mu)$ is always denoted by the symbol $\langle, \rangle$.

Let $\scrI_{1}(E,\mu)$ be the ideal of trace class operators
${\widetilde K}\colon L_2(E,\mu)\to L_2(E,\mu)$ (see volume~1 of~\cite{ReedSimon} for
the precise definition); the symbol
$||{\widetilde K}||_{\scrI_{1}}$ will stand for the
$\scrI_{1}$-norm of the operator ${\widetilde K}$.
Let $\scrI_{2}(E,\mu)$ be the
ideal of Hilbert-Schmidt operators ${\widetilde K}\colon
L_2(E,\mu)\to L_2(E,\mu)$; the symbol $||{\widetilde
K}||_{\scrI_{2}}$ will stand for the $\scrI_{2}$-norm of
the operator ${\widetilde K}$.

Let  $\mathscr I_{1,  \mathrm{loc}}(E,\mu)$ be the space of operators  $K\colon L_2(E,\mu)\to L_2(E,\mu)$
such that for any bounded Borel subset $B\subset E$
we have $$\chi_BK\chi_B\in\scrI_1(E,\mu).$$
Again, we endow the space $\scrI_{1,  \mathrm{loc}}(E,\mu)$
with a countable family of semi-norms
\begin{equation}
\label{btrcl}
||\chi_BK\chi_B||_{\scrI_1}
\end{equation}
where, as before, $B$ runs through an exhausting family $B_n$ of bounded sets.
A locally trace class operator $K$ admits a {\it kernel}, for which, slightly abusing notation, we use the same symbol $K$.

\subsection{Determinantal Point Processes}

A Borel probability measure $\mathbb{P}$ on
$\Conf(E)$ is called
\textit{determinantal} if there exists an operator
$K\in\scrI_{1,  \mathrm{loc}}(E,\mu)$ such that for any bounded measurable
function $g$, for which $g-1$ is supported in a bounded set $B$,
we have
\begin{equation}
\label{eq1}
\mathbb{E}_{\mathbb{P}}\Psi_g
=\det\biggl(1+(g-1)K\chi_{B}\biggr).
\end{equation}
Here and elsewhere in similar formulas, $1$ stands for the identity operator.
The Fredholm determinant in~\eqref{eq1} is well-defined since
$K\in \scrI_{1, \mathrm{loc}}(E,\mu)$.
The equation (\ref{eq1}) determines the measure $\Prob$ uniquely.
For any
pairwise disjoint bounded Borel sets $B_1,\dotsc,B_l\subset E$
and any  $z_1,\dotsc,z_l\in {\mathbb C}$ from (\ref{eq1}) we  have
$\mathbb{E}_{\mathbb{P}}z_1^{\#_{B_1}}\dotsb z_l^{\#_{B_l}}
=\det\biggl(1+\sum\limits_{j=1}^l(z_j-1)\chi_{B_j}K\chi_{\sqcup_i B_i}\biggr).$

For further results and background on determinantal point processes, see e.g.   \cite{Ghosh}, \cite{HoughEtAl}, \cite{Lyons},
 \cite{Lytvynov},  \cite{ShirTaka1}, \cite{ShirTaka2}, \cite{Soshnikov}.

If $K$ belongs to
$\scrI_{1, \text{loc}}(E,\mu)$, then, throughout the paper, we denote
the corresponding determinantal measure by
$\mathbb{P}_K$. Note that $\mathbb{P}_K$ is uniquely defined by~$K$,
but different operators may yield the same measure.
By a theorem due to Macch{\` \i} and Soshnikov ~\cite{Macchi}, \cite{Soshnikov} and Shirai-Takahashi \cite{ShirTaka0}, any
Hermitian positive contraction that belongs
to the class~$\scrI_{1, \text{loc}}(E,\mu)$ defines a determinantal point process.
For the purposes of this paper, we will only be interested in determinantal point processes given by operators of orthogonal projection; in the case of a discrete phase space, there is a standard procedure of doubling the phase space (see e.g. the Appendix in \cite{BOO}) that reduces the case of contractions to the case of projections.
\subsection{Weak convergence in the space of configurations.}
The space $\Conf(E)$ is  endowed with a natural structure of a complete
separable metric space, and the space   of finite Borel measures on the space of configurations is consequently also
a complete separable metric space with respect to the weak topology.  If $\varphi:E\to\mathbb{R}$ be a compactly supported continuous  function, then
 a measurable function
$\#_{\varphi}:\Conf(E)\to\mathbb{R}$ is introduced by the formula
$\#_{\varphi}(X)=\sum\limits_{x\in X}\varphi(x).$
Theorem 11.1.VII in vol.2 of \cite{DVJ} states that if $\Prob_n, n\in{\mathbb N}$ and $\Prob$ are Borel probability measures on $\Conf(E)$, then the measures $\Prob_n$ converge to $\Prob$ weakly as $n\to\infty$ if and only if for any finite collection $\varphi_1, \dots, \varphi_l$ of  continuous  functions  with disjoint compact supports
the joint distributions of the random variables
$\#_{\varphi_1}, \dots, \#_{\varphi_l}$ with respect to $\Prob_n$ converge, as $n\to\infty$,
to the joint distribution of  $\#_{\varphi_1}, \dots, \#_{\varphi_l}$ with respect to $\Prob$; convergence of joint distributions being understood according to the weak topology
on the space of Borel probability measures on ${\mathbb R}^l$. From the definition of determinantal point processes, we now have
\begin{proposition}\label{loctracetoweak}
Let $K_n$, $n\in {\mathbb N}$, $K$ be locally trace class operators acting in $L_2(E, \mu)$ and inducing determinantal mesures $\Prob_{K_n}$,  $\Prob_K$ .
If $K_n\to K$ in  $\scrI_{1, \text{loc}}(E,\mu)$ as $n\to\infty$, then $\Prob_{K_n}\to \Prob_K$
weakly in the space of probability measures on $\Conf(E)$.
\end{proposition}

\subsection{The product of a determinantal measure and a multiplicative functional.}
We start by recalling the results of \cite{Buf-CIRM}, \cite{Buf-umn}   showing that the product  of a determinantal measure with a multiplicative functional is, after normalization,  again a determinantal measure, whose kernel is found explicitly.

Let $g$ be a non-negative measurable function on $E$. If the operator $1+(g-1)K$ is invertible, then we set
$${\mathfrak B}(g, K)=g K(1+{(g-1)}K)^{-1},\qquad
{\widetilde {\mathfrak B}}(g, K)= {\sqrt{g}}K(1+{(g-1)}K)^{-1}{\sqrt{g}}.$$
By definition, ${\mathfrak B}(g,K),{\widetilde {\mathfrak B}}(g,K)\in \scrI_{1,\text{loc}}(E,\mu)$ since $K\in \scrI_{1,\text{loc}}(E,\mu)$, and, if $K$ is self-adjoint, then
so is ${\widetilde {\mathfrak B}}(g,K)$ (this follows from the elementary identity $p(1-qp)^{-1}=(1-pq)^{-1}p$ that holds for arbitrary bounded operators $p,q$ provided both sides are well-defined).

We now quote Proposition 2.1 in \cite{Buf-CIRM}.
\begin{proposition}
\label{pr1-bis}
Let $K\in \scrI_{1,\mathrm{loc}}(E,\mu)$ be a self-adjoint positive contraction, and let $\mathbb{P}_K$ be
the corresponding determinantal measure on $\Conf(E)$. Let
$g$ be a nonnegative bounded measurable
function on~$E$ such that
\begin{equation}
\label{gkint}
\sqrt{g-1}K\sqrt{g-1}\in \scrI_{1}(E,\mu),
\end{equation}
and that the operator
$1+{(g-1)}K$ is invertible.  Then the operators
${\mathfrak B}(g, K), {\widetilde {\mathfrak B}}(g, K)$ induce on $\Conf(E)$ a determinantal measure
$\mathbb{P}_{{\mathfrak B}(g, K)}=\mathbb{P}_{{\widetilde {\mathfrak B}}(g, K)}$ satisfying
\begin{equation}\label{mult1}
\mathbb{P}_{{\mathfrak B}(g, K)}=\displaystyle \frac{\displaystyle \Psi_g\mathbb{P}_K}
{\displaystyle \displaystyle \int\limits_{\Conf(E)}\Psi_g\,d\mathbb{P}_K}.
\end{equation}

\end{proposition}
{\bf Remark.} Here and elsewhere, we write $\sqrt{g-1}K\sqrt{g-1}$ instead of the longer 
and more formal $\mathrm{sgn} (g-1)\sqrt{|g-1|}K\sqrt{|g-1|}$.

{\bf Remark.} Of course, from (\ref{gkint}) and the invertibility of the operator $1+{(g-1)}K$ we have
 $\Psi_g\in L_1(\Conf(E),\mathbb{P}_K)$ and
$$
\displaystyle \int\Psi_g\,d\mathbb{P}_K=\det\Bigl(1+\sqrt{g-1}K\sqrt{g-1}\Bigr)>0,
$$
so the right-hand side of (\ref{mult1}) is well-defined.

For the reader's convenience, we recall the proof of Proposition \ref{pr1-bis}  in the case when the assumption (\ref{gkint}) is replaced (cf. \cite {Buf-umn}) by a simpler assumption $$(g-1)K\in \scrI_{1}(E,\mu);$$
 for  the general case, see  the proof of Proposition 2.1 in \cite{Buf-CIRM}.
Take a bounded measurable function $f$ on $E$ such that
$(f-1)K\in\scrI_{1}(E,\mu)$; for example, one may take $f$ that is different from $1$ on a bounded set.
We have $(fg-1)K\in\scrI_{1}(E,\mu)$ since  $(f-1)K\in \scrI_{1}(E,\mu)$, $(g-1)K\in \scrI_{1}(E,\mu)$.
By definition, we  have
\begin{multline}\label{mult-det1}
{\mathbb E}_{{\mathbb P}_K} \Psi_f\Psi_g=\det(1+(fg-1)K)=\\=\det (1+(f-1)gK(1+(g-1)K)^{-1})\det (1+(g-1)K).
\end{multline}
We rewrite (\ref{mult-det1}) in the form
$$
\displaystyle \frac{{\mathbb E}_{{\mathbb P}_K} \Psi_f\Psi_g}{{\mathbb E}_{{\mathbb P}_K}\Psi_g}
=\det(1+(f-1){\mathfrak B}(g,K))=
\det(1+{(f-1)}{\widetilde {\mathfrak B}}(g,K)).
$$
Since a probability measure on the space of configurations is uniquely determined by the values of multiplicative functionals corresponding to all bounded functions $f$ that are different from $1$
on a bounded set, formula (\ref{mult-det1}) implies Proposition \ref{pr1-bis}.

\subsection{Projections and subspaces.}

Let $L \subset L_2(E, \mu)$ be a closed subspace, let $\Pi$ be the corresponding projection operator, assumed to be locally of trace class, and let $\mathbb{P}_{\Pi}$ the corresponding determinantal measure.
Our aim is to determine how the measure $\mathbb{P}_{\Pi}$ changes if the subspace $L$ is multiplied by a function.
We start with the following clear
\begin{proposition}
Let $\alpha(x)$ be a measurable function such that $|\alpha(x)|=1$ $\mu$-almost surely.
Then the operator of orthogonal projection onto the subspace $\alpha(x)L$ induces the same determinantal measure $\mathbb{P}_{\Pi}$.
\end{proposition}
Proof.
Indeed, if $\Pi(x,y)$ is the kernel of the operator $\Pi$, then the kernel of the new operator has the form

$$
\displaystyle \frac{\alpha(x) \Pi(x,y)}{\alpha(y)},
$$
and such gauge transformations do not change the determinantal measure.

\begin{proposition}
Let $g$ be a non-negative bounded function on $E$ such that the operator $1+(g-1)\Pi$ is invertible.
Then the operator
\begin{equation}\label{pig-def}
\Pi^g=\sqrt{g}\Pi(1+(g-1)\Pi)^{-1} \sqrt{g}
\end{equation}
is the operator of orthogonal projection onto the closure of the subspace $\sqrt{g}L$.
\end{proposition}
\begin{proof}
 First, let $\widetilde{\varphi} \in \sqrt{g}L$, that is, $\widetilde{\varphi} = \sqrt{g}\varphi, \varphi \in L$. Since $\varphi \in L$, we have
$$ (1+(g-1)\Pi)\varphi = g\varphi, $$
whence
$$  (1+(g-1)\Pi)^{-1} \sqrt{g}\widetilde{\varphi} = \varphi, $$
and finally  $$\Pi^{g}\widetilde{\varphi} = \widetilde{\varphi}$$ as desired. \par
Now take $\varphi$ to be orthogonal to the subspace $\sqrt{g}L$. Since $g$ is real-valued, we have $\sqrt{g}\varphi \in L^{\bot}$, whence
$(1+(g-1)\Pi)\varphi= \sqrt{g}\varphi$ and, consequently,
 $(1+(g-1)\Pi)^{-1} \sqrt{g}\varphi = \varphi$. We thus have $\Pi^g \varphi = 0$, and the proposition is proved completely.
\end{proof}

We prepare a useful proposition for proving local trace-class convergence of sequences of
 operators $\Pi^g$.

\begin{proposition}\label{pign-conv}
Let $g_n$, $n\in {\mathbb N}$ be a sequence of nonnegative Borel functions satisfying
\begin{equation}\label{bddgn}
\sup\limits_{n\in {\mathbb N},x\in E} g_n(x)<+\infty,
\end{equation}
 and assume that
\begin{enumerate}
\item as $n\to\infty$ the sequence $g_n$ converges, $\mu$-almost surely, to a Borel function $g$ ;
\item  the operator $1+(g-1)\Pi$ is invertible ;
\item  $||(g_n-g)\Pi||\to 0$.
\end{enumerate}
Then
$$
\Pi^{g_n}\to \Pi^g \ \mathrm{in}  \scrI_{1, \mathrm{loc}}(E,\mu).
$$
\end{proposition}
\begin{proof}
Let $B$ be a bounded set and write $$\chi_B\Pi^g\chi_B=(\sqrt g \chi_B\Pi) (\Pi(1+(g-1)\Pi)^{-1}\Pi)(\Pi\chi_B\sqrt{g}).$$
If $A$ ranges in the space of Hilbert-Schmidt operators and $D$ in the space of bounded self-adjoint operators endowed with the operator norm, then
the correspondence $(A, D)\to A^*DA$ induces a continuous mapping into the space
of trace-class operators.  By our third assumption, the operators $1+(g_n-1)\Pi$ converge to
$1+(g-1)\Pi$ in the norm topology, consequently,  for large $n$ the operators $1+(g_n-1)\Pi$ are invertible, the inverses also converge in norm, and, finally,
$\Pi(1+(g_n-1)\Pi)^{-1}\Pi\to \Pi(1+(g-1)\Pi)^{-1}\Pi$ in the norm topology. Uniform boundedness (\ref{bddgn})
of $g_n$ and pointwise convergence of $g_n$ imply that $\sqrt g_n \chi_B\Pi\to \sqrt g \chi_B\Pi$
in the Hilbert-Schmidt norm, and the proof is complete.
\end{proof}

\subsection{Normalized multiplicative functionals}

If the multiplicative functional $\Psi_g$ is $\Prob_{\Pi}$-integrable, then we introduce the {\it normalized} multiplicative functional
${\overline \Psi}_g$ by the formula
\begin{equation}\label{psig-def}
{\overline \Psi}_g=\displaystyle \frac{\Psi_g}{\displaystyle \int\limits_{\Conf(E)} \Psi_gd\Prob_{\Pi}}.
\end{equation}

We reformulate Proposition 2.1 in \cite{Buf-CIRM} in our  new notation (\ref{pig-def}), (\ref{psig-def}):
\begin{proposition}\label{mult-prelim}
If $g$ is a bounded Borel function on $E$ such that $$\sqrt{g-1}\Pi\sqrt{g-1}\in \scrI_{1}(E,\mu)$$
and the operator $1+(g-1)\Pi$ is invertible, then  the subspace $\sqrt{g}L$ is closed, the normalized multiplicative functional
${\overline \Psi}_g$ is well-defined,   and we have
\begin{equation}\label{pig-psig}
{\overline \Psi}_g\Prob_{\Pi}=\Prob_{\Pi^g}.
\end{equation}
\end{proposition}

Note that closedness of the subspace $\sqrt{g}L$ is immediate from the invertibility of the operator $1+(g-1)\Pi$: indeed, the operator $1+(g-1)\Pi$ takes the subspace $L$ to the subspace $gL$, which is consequently closed.
 Since the function $g$ is bounded from above, the subspace $\sqrt{g}L$ is, a fortiori, closed as well: indeed, if $\psi_n\in L$ are such that
$\sqrt{g}\psi_n\to \psi$, then, by boundedness of $g$, we also have $g\psi_n\to \sqrt{g}\psi$, so, by closedness
of the space $gL$, we have $\sqrt{g}\psi\in gL$ whence also $\psi\in \sqrt{g}L$ (this implication is valid  even if the set $\{x:g(x)=0\}$ has positive measure since, by definition, $\psi=\chi_{\{x:g(x)>0\}}g$).

A key point in the argument of this paper is that the normalized multiplicative functional (\ref{psig-def}) can be defined, in such a way that the formula (\ref{pig-psig}) still holds,  even when the multiplicative functional $\Psi_g$ itself is not defined, see Proposition \ref{multreg-main} below.

\subsection{On the subspace $\sqrt{g}L$.}
We now give another sufficient condition for closedness of the subspace $\sqrt{g}L$ for a bounded function $g$ .
If the function $g$ is bounded away from zero in addition to being bounded from above, then the subspace $\sqrt{g}L$ is automatically closed. If  $\inf\limits_{x\in E} g=0$, then we prepare the following simple proposition.
Recall that $||A||$ always stands for the usual operator norm of a bounded  operator $A$.
\begin{proposition}
Let $\Pi$ be an operator of orthogonal projection onto a closed subspace $L\subset L_2(E, \mu)$,
and let $C\subset E$ be a Borel subset such that $||\chi_{E\setminus C}\Pi||<1$. Then the subspace $\chi_CL$ is closed and  the natural  restriction map $\varphi \to \chi_C\varphi$ induces an isomorphism of Hilbert spaces $L$ and $\chi_CL$.
\end{proposition}
\begin{proof}
Indeed, our assumptions imply the existence of a positive constant $\alpha$ such that for any
$\varphi\in L$ we have
\begin{equation}\label{alpha-est-l}
||\chi_{E\setminus C}\varphi||\leq \alpha||\chi_{C}\varphi||,
\end{equation}
and the proposition follows.
\end{proof}
\begin{corollary}\label{gclose-cor}
Let $g$ be a bounded nonnegative Borel function on $E$ and $C\subset E$ a Borel subset such that
\begin{enumerate}
\item $||\chi_{E\setminus C}\Pi||<1$;
\item the function $g$ in restriction to $C$ is bounded away from $0$.
\end{enumerate}
Then the subspace $\sqrt{g}L$ is closed.
\end{corollary}
\begin{proof}
If $g|_C>\varepsilon>0$ and $g<M<+\infty$ on the whole space $E$, then the inequality (\ref{alpha-est-l})
implies, for any $\psi\in \sqrt{g}L$, the estimate
\begin{equation}\label{alpha-est-l-2}
||\chi_{E\setminus C}\psi||\leq M\varepsilon^{-1}\alpha||\chi_{C}\psi||.
\end{equation}
The subspace $\chi_CL$ is closed by our second assumption, and so is $\chi_C\sqrt{g}L$. By
(\ref{alpha-est-l-2}),  for all $\psi\in \sqrt{g}L$, the natural restriction map  $\psi\to \chi_C\psi$
is invertible with bounded inverse, and the proposition follows.
\end{proof}

\begin{proposition}\label{gless1bd}
Let $g:E\to [0,1]$ be a Borel function such that
\begin{enumerate}
\item
$\tr \left(\chi_{\{x\in E: g(x)<1\}} \Pi \chi_{\{x\in E: g(x)<1\}} \right)<+\infty$,
\item a function $\varphi\in L$ satisfying $\chi_{\{x: g(x)<1\}}\varphi=\varphi$ must be the zero function.
\end{enumerate}
Then all the conclusions of Proposition \ref{mult-prelim} hold for the function $g$.
\end{proposition}
\begin{proof}
Our assumptions imply
$||\chi_{\{x\in E: g(x)<1\}} \Pi||<1$,   whence the operator $1+(g-1)\Pi$  is invertible.
The space $\sqrt{g}L$ is closed by Corollary \ref{gclose-cor}.
\end{proof}

\subsection{On invertibility of the operator $1+(g-1)\Pi$.}

Take  a nonnegative bounded Borel function $g$ such that the operator $\sqrt{g-1}\Pi$ is Hilbert-Schmidt. It follows that the operator $(g-1)\Pi$ is also Hilbert-Schmidt, while the operator
$\sqrt{g-1}\Pi\sqrt{g-1}$ is trace-class. In particular, since all these operators are compact, invertibility
of the operator $1+(g-1)\Pi$  is equivalent to its injectivity, and the same is true for $1+\sqrt{g-1}\Pi\sqrt{g-1}$.
We start with a simple remark.
\begin{proposition}
\begin{enumerate}
\item The  invertibility of the operator $1+(g-1)\Pi$  is equivalent to the invertibility of the operator $1+\sqrt{g-1}\Pi\sqrt{g-1}$.
\item If $\sup|g-1|<1$, then the operator $1+(g-1)\Pi$ is invertible.
\end{enumerate}
\end{proposition}
\begin{proof}
Indeed, if $\sqrt{g-1}\Pi\sqrt{g-1}\varphi=-\varphi$, then $(g-1)\Pi(\sqrt{g-1}\varphi)=-\sqrt{g-1}\varphi$.
Conversely, if $(g-1)\Pi\varphi=-\varphi$, then the function $\psi=\sqrt{g-1}\Pi\varphi$ satisfies
$\sqrt{g-1}\psi=-\varphi$ and, consequently, $\sqrt{g-1}\Pi\sqrt{g-1}\psi=-\psi$, and the first item is proved.  The second item is clear since  if $\sup|g-1|<1$, then  $||(g-1)\Pi||<1$.
\end{proof}

We next show that perturbing a positive function $g$ on a bounded set does not change the  invertibility of the operator $1+(g-1)\Pi$.

\begin{proposition}\label{bdd-pert}
Let $g_1, g_2$ be positive bounded Borel functions  such that the operators $\sqrt{g_1-1}\Pi$, $\sqrt{g_2-1}\Pi$ are Hilbert-Schmidt and the set $\{x: g_1(x)\neq g_2(x)\}$ is bounded.
Then the  invertibility of the operator $1+(g_1-1)\Pi$  is equivalent to the invertibility of the operator
 $1+(g_2-1)\Pi$.
\end{proposition}
\begin{proof}
We have
$
\ee_{\Prob_{\Pi}} \Psi_{g_i}=\det(1+\sqrt{g_i-1}\Pi\sqrt{g_i-1}), \ i=1,2.
$
Since our functions are positive,  the operator $1+(g_i-1)\Pi$ is invertible if and only if  the multiplicative functional $\Psi_{g_i}$is positive  on a subset of positive $\Prob_{\Pi}$-measure of $\Conf(E)$,
or, in other words, if and only if the   infinite product $\prod\limits_{x\in X} g_i(x)$ converges absolutely
 to a   positive limit  with positive $\Prob_{\Pi}$-probability.  Since both  functions are positive and the set $\{x: g_1(x)\neq g_2(x)\}$ is bounded, for $\Prob_{\Pi}$-almost any configuration $X$ the functions $g_1$ and $g_2$ are equal except on finitely many particles of $X$,  and   the   infinite products $\prod\limits_{x\in X} g_i(x)$, $i=1,2$, either both converge or both diverge.
\end{proof}
\begin{proposition}\label{restr-inv}
If $g$ is a nonnegative bounded Borel function satisfying
\begin{enumerate}
\item  the operator  $(g-1)\Pi$ is compact;
\item the operator $1+(g-1)\Pi$ is invertible,
\end{enumerate}
 then, for any subset $C\in E$, the function $g^C=\chi_{E\setminus C}+(g-1)\chi_C$ also satisfies the same conditions 1,2.
\end{proposition}
\begin{proof}
Compactness of $(g^C-1)\Pi=\chi_C(g-1)\Pi$ is clear. If $1+(g^C-1)\Pi$ fails to be invertible and has a nontrivial kernel,
then there exists a function $\psi$ satisfying $\psi+\chi_C(g-1)\Pi\psi=0$. It follows that $\psi=\chi_C\psi$ and that
$\psi+(g-1)\Pi\psi=0$, a contradiction that completes the proof.
\end{proof}

\subsection{Inverting self-adjoint operators}

It is sometimes more convenient to  take the inverse of $1+\sqrt{g-1}\Pi\sqrt{g-1}$ rather than that of $1+(g-1)\Pi$, and we
rewrite the formula (\ref{pig-def}) in the following way.
Let $\alpha>0$, let $g$ be a Borel function, perhaps unbounded from above,  on $E$, satisfying $g>1+\alpha$ and such that the operator
$\sqrt{g-1}\Pi$ is Hilbert-Schmidt. Then $\sqrt{g}L$ is a subspace of $L_2(E, \mu)$, automatically a closed one; we let $\Pi^g$ be the operator of   orthogonal projection onto $\sqrt{g}L$ and set $\Pi_g=\sqrt{g-1}\Pi\sqrt{g-1}$.
The operator $1+\Pi_g$ is then automatically invertible.
\begin{proposition}\label{pigprop}
We have
\begin{equation}\label{pig-pig}
\Pi^g=\sqrt{1+(g-1)^{-1}}\Pi_g(1+\Pi_g)^{-1}\sqrt{1+(g-1)^{-1}}.
\end{equation}
\end{proposition}
\begin{proof}
If $g$ is bounded, then, noting that invertibility of $1+\Pi_g$ implies that of $1+(g-1)\Pi$,  we again use the identity $p(1+qp)^{-1}=(1+pq)^{-1}q$, valid for arbitrary bounded operators $p,q$ once both sides are well-defined.  If $g$ is unbounded, then we approximate $g$ from below by a sequence of bounded functions and pass to the limit in the space of Hilbert-Schmidt operators.
\end{proof}
{\bf Remark.} While this proposition is sufficient for our purposes, the equality \eqref{pig-pig} still holds under much milder assumptions on $g$; the key point is to
ensure invertiblity of the operators $1+\Pi_g$, $1+(g-1)\Pi$.
The formula (\ref{pig-pig}) implies the following analogue of Proposition \ref{pign-conv}.
\begin{proposition}\label{pign-sa-conv}
Let $\alpha>0$, and let $g$ be a Borel fucntion on $E$ satisfying $g>1+\alpha$, $\sqrt{g-1}\Pi\in \scrI_2(E, \mu)$.
Let $g_n$, $n\in {\mathbb N}$ be a sequence of nonnegative bounded Borel functions satisfying
$
1+\alpha\leq \inf\limits_{n\in {\mathbb N},x\in E} g_n(x),
$
 and assume that
 $\sqrt{g_n-1}\Pi\to \sqrt{g-1}\Pi\in \scrI_2(E, \mu)$.
Then
$$
\Pi^{g_n}\to \Pi^g \ \mathrm{in}  \scrI_{1, \mathrm{loc}}(E,\mu).
$$
\end{proposition}
\begin{proof} The Hilbert-Schmidt convergence of the sequence $\sqrt{g_n-1}\Pi$ implies the trace class, and, a fortiori, norm,  convergence of the sequence $\Pi_{g_n}$. The proof is now concluded in the same way as that of Proposition \ref{pign-conv}.
\end{proof}

\subsection{Palm Measures of Determinantal Point Processes.}
Palm measures of determinantal point processes admit the following characterization.
As above, let $\Pi\in \scrI_{1,\mathrm{loc}}(E,\mu)$ be the operator of orthogonal projection onto a closed subspace $L\subset L_2(E, \mu)$.
For $q\in E$ satisfying $\Pi(q,q)\neq 0$, introduce a kernel $\Pi^q$ by the formula
\begin{equation}\label{piqdef}
\Pi^q(x,y)=\Pi(x,y)-\displaystyle \frac{\Pi(x,q)\Pi(q,y)}{\Pi(q,q)}.
\end{equation}
If $\Pi(q,q)=0$, then we also have $\Pi(x,q)=\Pi(q,y)=0$ almost surely with respect to $\mu$, and we set $\Pi^q=\Pi$.

The operator $\Pi^q$ defines an orthogonal projection onto the subspace
$$
L(q)=\{\varphi\in L: \varphi(q)=0\}
$$
of functions in $L$ that assume the value zero at the point $q$; the space $L(q)$ is well-defined by Assumption \ref{ass-ker}; in other words, $L(q)$ is the orthogonal complement of $v_q$ in $L$.
Iterating, let $q_1, \ldots, q_l\in E$ be distinct and set
$$
L(q_1, \dots, q_l)=\{\varphi\in L: \varphi(q_1)=\ldots=\varphi(q_l)=0\},
$$
and let $\Pi^{q_1, \dots, q_l}$ be the operator of orthogonal projection onto the subspace $L(q_1, \dots, q_l)$.
Shirai and Takahashi \cite{ShirTaka1} have proved
\begin{proposition}[Shirai and Takahashi \cite{ShirTaka1}]\label{palm-char}
For  any $l\in {\mathbb N}$ and for $\rho_l$-almost every $l$-tuple $q_1, \dots, q_l$ of distinct points in $E$,
the iterated reduced Palm measure $\Prob_{\Pi}^{q_1, \dots, q_l}$ is given by the formula
\begin{equation}\label{ppiq-it}
\Prob_\Pi^{q_1, \dots, q_l}= \Prob_{\Pi^{q_1, \dots, q_l}}.
\end{equation}
\end{proposition}

{\bf Remark.} Shirai and Takahashi \cite{ShirTaka1} have in fact established the
formula (\ref{ppiq-it}) for arbitrary positive self-adjoint
locally trace-class contractions; the formula (\ref{piqdef}) for the kernel stays the same. Note that in the discrete case the formula for contractions is a corollary of the formula for projection operators, since formula (\ref{ppiq-it}) is local (both in the discrete and the continuous cases) and contractions are reduced to projections by  doubling the phase space (see e.g. the Appendix in Borodin-Okounkov-Olshanski \cite{BOO}).

\subsection{Conditional measures in the discrete case}

In this subsection, we consider the discrete case, in which the space $E$ is a countable set endowed with the discrete topology, and the measure $\mu$ is the counting measure. In this case, the reduced Palm measure
$\Prob^{ q}$
of a point process $\Prob$ on $\Conf(E)$
can be described as follows: one takes the conditional measure of $\Prob$ on the subset of configurations containing a particle at position $q$, and then one removes the particle at $q$; more formally, $\Prob^{ q}$ is the push-forward of the said conditional measure under the operation that to a configuration $X$ containing the particle at $q$ assigns the configuration $X\setminus \{q\}$.

In the discrete case we also have a dual construction: let $\Prob^{\breve q}$ be the conditional measure of $\Prob$ with respect to the event that there is no particle at position $q$. More formally,
set
$$
\Conf(E; E\setminus \{q\})=\{X\in \Conf(E): q\notin X\},
$$
and write
$$
\Prob^{\breve q}=\displaystyle \frac{\Prob|_{\Conf(E; E\setminus \{q\})}}{\Prob(\Conf(E; E\setminus \{q\}))}
$$
be the normalized restriction of $\Prob$ onto the subset $\Conf(E; E\setminus \{q\})$.

We have a dual to Proposition \ref{palm-char}.

\begin{proposition} \label{palm-hole}
Let $q \in E$ be such that $\mu(\{q\})> 0$.
Then the operator of orthogonal projection onto the subspace $\chi_{E \setminus q}L$ has the kernel $\Pi^{\breve{q}}$ given by the formula
\begin{equation}\label{pqbreve}
\Pi^{\breve{q}}(x, y) = \Pi(x,y) + \displaystyle \frac{\Pi(x, q)\Pi(q, y)}{1-\Pi(q,q)}, x \neq q, y \neq q;\\
\end{equation}
\begin{equation}
\Pi^{\breve{q}}(x, q) = \Pi^{\breve{q}}(q, y) = 0, x,y\in E.
\end{equation}
\end{proposition}
Proof. This is a particular case of Corollary 6.4 in Lyons \cite{Lyons}; see also Shirai-Takahashi ~\cite{ShirTaka1}, ~\cite{ShirTaka2}.

Given $l\in {\mathbb N}$, $m<l$ and an $l$-tuple $(p_1, \dots, p_l)$, of distinct points in $E$, recall that we have introduced a subspace  $L(p_1, \dots, p_m, \breve {p}_{m+1}, \dots, \breve {p}_l)$ by the formula
\begin{equation}\label{lppbreve-1}
  L(p_1, \dots, p_m, \breve {p}_{m+1}, \dots, \breve {p}_l) =
\{  \chi_{E\setminus \{{p}_{m+1}, \dots, {p}_l\}}\varphi : \varphi \in L, \varphi (p_1) =\dots=\varphi(p_m)= 0 \}.
\end{equation}

Let $\Pi^{p_1, \dots, p_m, \breve {p}_{m+1}, \dots, \breve {p}_l}$ be the operator of orthogonal projection onto the
subspace $L(p_1, \dots, p_m, \breve {p}_{m+1}, \dots, \breve {p}_l)$. The corresponding determinantal measure
$\Prob_{\Pi^{p_1, \dots, p_m, \breve {p}_{m+1}, \dots, \breve {p}_l}}$ admits the following characterization.
Recall that $${\bf C}(p_1, \dots, p_m, \breve{p}_{m+1}, \dots, \breve{p}_l)$$ is the set of configurations on $E$ containing
exactly one particle in each of the positions $p_1, \dots, p_m$ and no particles in the positions $p_{m+1}, \dots, p_l$.
There is a natural  {\it erasing bijection} between ${\bf C}(p_1, \dots, p_m, \breve{p}_{m+1}, \dots, \breve{p}_l)$ and
${\bf C}({\breve p}_1, \dots, {\breve p}_m, \breve{p}_{m+1}, \dots, \breve{p}_l)$ obtained by erasing the particles
in positions $p_1, \dots, p_m$.

\begin{proposition}
Consider the normalized restriction of $\Prob_{\Pi}$  onto the set ${\bf C}(p_1, \dots, p_m, \breve{p}_{m+1}, \dots, \breve{p}_l)$. The push-forward of this normalized restriction onto the set ${\bf C}({\breve p}_1, \dots, {\breve p}_m, \breve{p}_{m+1}, \dots, \breve{p}_l)$ under the erasing bijection is the measure $\Prob_{\Pi^{p_1, \dots, p_m, \breve{p}_{m+1}, \dots, \breve{p}_l}}$.
\end{proposition}
Proof. Again, this is a reformulation of Corollary 6.4 in Lyons \cite{Lyons}; see also Shirai-Takahashi ~\cite{ShirTaka1}, ~\cite{ShirTaka2}.

\subsection{Action of Borel automorphisms on point processes.}

Let $T: E \to E$ be an invertible measurable map such that for any bounded set $B \subset E$ the set
$T(B)$ is also bounded. The map $T$ naturally acts on the space of configurations $\mathrm{Conf}(E)$: namely, given $X \in \mathrm{Conf}(E)$ we set
$$
	T(X) =  \{Tx, x \in X\}
$$
Note that, by our assumptions, $T(X)$ is a well-defined configuration on $E$; slightly abusing notation, we
thus keep the same symbol $T$ for the induced action on the space of configurations.

Let $\mathbb{P}$ be a probability measure on $\mathrm{Conf}(E)$. We assume that $\mathbb{P}$ admits correlation measures of all orders, and, for $l \in \mathbb{N}$, we let $\rho_l$ be the $l$-th correlation measure of the point process $\mathbb{P}$. The $l$-th Cartesian power of $T$ naturally acts on the measure $\rho_l$, and, slightly abusing notation, we denote the resulting measure by $\rho_l\circ T$. The measure $\rho_l\circ T$ is, of course, the $l$-th correlation measure of the point process $\mathbb{P}\circ T$, the push-forward of the measure $\mathbb{P}$ under the induced action of the automorphism $T$ on the space of configurations.

We now prove a simple general statement: if  for  a point process $\Prob$ and an arbitrary fixed  $l\in {\mathbb N}$, the reduced Palm measures
 corresponding to different $l$-tuples of points are equivalent, then for any Borel automorphism $T$
acting by the identity beyond a bounded set, the measures
 $\mathbb{P}$ and $\mathbb{P}\circ T$ are also equivalent, and the Radon-Nikodym derivative is found explicitly in terms of the Radon-Nikodym derivatives of the reduced Palm measures. More precisely, we have the following
\begin{proposition}\label{RN-T}
Let $T:E \to E$ be a Borel automorphism admitting a bounded subset $B \subset E$ such that $T(x) = x$ for all $x \in E\backslash B$. Assume that
\begin{enumerate}
\item for any $l \in \mathbb{N}$, the correlation measures $\rho_l$ and $\rho_l\circ T$
are equivalent;
\item for any two collections $\{q_1,...,q_l\}$ and $\{q'_1,...,q'_l\}$ of distinct points of $E$, the measures ${{\mathbb{P}}}^{q_1,...,q_l}$ and ${{\mathbb{P}}}^{q'_1,...,q'_l}$ are equivalent.
\end{enumerate}
Then the measures $\mathbb{P}$ and $\mathbb{P}\circ T$ on $\mathrm{Conf}(E)$ are equivalent, and for $\mathbb{P}$-almost  every configuration $X \in \mathrm{Conf}(E)$ such that $X \cap B=\{q_1,...,q_l\} $ we have
$$
	\displaystyle \frac{d\mathbb{P}\circ T}{d\mathbb{P}}(X) = \displaystyle \frac{d{{\mathbb{P}}}^{T{q_1},...,T{q_l}}}{d{{\mathbb{P}}}^{q_1,...,q_l}}
	(X\setminus \{q_1, \dots, q_l\}) \times \displaystyle \frac{d\rho_l\circ T}{d\rho_l}(q_1,...,q_l).
$$
\end{proposition}
\begin{proof}
Let $l\in \mathbb{N}$, let $\mathbb{P}$, $\widetilde{\mathbb{P}}$ be
 probability measures on $\mathrm{Conf}(E)$ admitting correlation measures of order $l$, denoted, respectively, $\rho_l$,  $\widetilde {\rho}_l$ . Let $\Prob^{q_1, \dots, q_l}$, ${\widetilde {\Prob}}^{q_1, \dots, q_l}$ stand for the
respective reduced Palm measures. The symbol $\ll$ denotes absolute continuity
of measures. Proposition \ref{RN-T}  follows from

\begin{proposition}\label{pp-cont}
If $\widetilde {\rho}_l\ll\rho_l$ and ${\widetilde {\Prob}}^{q_1, \dots, q_l}\ll\Prob^{q_1, \dots, q_l}$ for $\rho_l$-almost any distinct $q_1, \dots, q_l\in E$, then also $\widetilde{\mathbb{P}}\ll\mathbb{P}$ and for $\Prob$-almost any $X\in\Conf(E)$ and any
$l$ particles $q_1, \dots, q_l\in X$ we have

$$
	\displaystyle \frac{d\widetilde{\mathbb{P}}}{d\mathbb{P}} = \displaystyle \frac{d{\widetilde{\mathbb{P}}}^{{q_1},...,{q_l}}}{d{{\mathbb{P}}}^{q_1,...,q_l}}
	(X\setminus \{q_1, \dots, q_l\}) \times \displaystyle \frac{d\widetilde {\rho}_l}{d\rho_l}(q_1,...,q_l).
$$
\end{proposition}
\begin{proof}
We use  Proposition \ref{palm-cyl}.
Take $X_0\in\Conf(E)$ and $q_1^{(0)}, \dots, q_l^{(0)}\in X_0$; take disjoint bounded open sets $B^{(1)},...,B^{(l)}\subset E$, set $B=\cup B^{(i)}$ and
take an open set $U\subset E$
disjoint from all $B^{(i)}$ in such a way that
$q_i(0)\in B^{(i)}$ for all $i=1, \dots, l$ and $X_0\setminus \{q_1^{(0)}, \dots, q_l^{(0)}\}\subset U$. Let $\mathscr W$ be a neighbourhood of $X_0\setminus \{q_1^{(0)}, \dots, q_l^{(0)}\}$ in $\Conf(E)$
satisfying $\mathscr W\subset \{X\in\Conf(E): X\subset U\}$.
Introduce a neighbourhood $\mathscr Z$ of $X_0$ by setting

\begin{equation}\label{bngen}
\mathscr Z = \{ X \in\Conf(E): {\#}_{B^{(1)}}(X)=\dots={\#}_{B^{(l)}}(X)=1, X|_{E\setminus B}\subset  {\mathscr W}\}.
\end{equation}
Sets given by (\ref{bngen}) form a basis of neighbourhoods of $X_0$.
We prepare a simple general lemma.
\begin{lemma}\label{lem-ac}
Let $\nu_1$ and $\nu_2$ be two Borel probability measures on a complete separable metric space $Y$.
Let $\mathscr V=\{V\}$ be a basis of neighbourhoods on $Y$.
Let $\varphi$ be a nonnegative function on $Y$ such that for any neighbourhood $V\in {\mathscr V}$ we have
\begin{equation}\label{nu12}
\nu_1(V)=\int\limits_V \varphi d\nu_2.
\end{equation}
Then $\nu_1\ll\nu_2$ and $d\nu_1/d\nu_2=\varphi$ almost surely with respect to $\nu_2$.
\end{lemma}
\begin{proof}
Indeed, if the equality (\ref{nu12}) holds for a basis of neighbourhoods, then, due to tightness,  it  holds for all open sets and then also
for all Borel sets.
\end{proof}

By definition of Palm measures and Proposition \ref{palm-cyl}, we have
$$
	\mathbb{P}(\mathscr Z) = \displaystyle \int\limits_{B^{(1)} \times ... \times B^{(l)}} {{\mathbb{P}}}^{q_1,...,q_l}(\mathscr W)
	d\rho_l(q_1,...,q_l).
$$

A similar formula holds for $\widetilde \Prob$.
Taking $X\in {\mathscr Z}$, setting $q_i=X\cap B^{(i)}$,  and using the assumption
${\widetilde {\Prob}}^{q_1, \dots, q_l}\ll\Prob^{q_1, \dots, q_l}$, we therefore have
$$
	{\widetilde \Prob}(\mathscr Z) = \displaystyle \int\limits_{\mathscr Z} \displaystyle \frac{d{\widetilde{\mathbb{P}}}^{{q_1},...,{q_l}}}{d{{\mathbb{P}}}^{q_1,...,q_l}}
	(X\setminus \{q_1, \dots, q_l\}) \times \displaystyle \frac{d\widetilde {\rho}_l}{d\rho_l}(q_1,...,q_l) d\Prob(X),
$$
Since the formula holds for an arbitrary neighbourhood of the form (\ref{bngen}),
 Proposition
\ref{pp-cont} follows from Lemma \ref{lem-ac}.
\end{proof}

We now derive Proposition \ref{RN-T} from
Proposition \ref{pp-cont}. As before, let $\Conf(E; E\setminus B)$ be the subset of those configurations on $E$ all whose particles lie in  $E\setminus B$.
Since the automorphism $T$ acts by the identity on $E\setminus B$, all configurations in the set $\Conf(E; E\setminus B)$ are fixed by $T$, and
we have
$$
\Prob^{q_1, \dots q_l}|_{\Conf(E; E\setminus B)}\circ T=\Prob^{Tq_1, \dots Tq_l}|_{\Conf(E; E\setminus B)}.
$$
By definition, we have
$X\setminus \{q_1, \dots, q_l\}\in \Conf(E; E\setminus B),$  whence, again using neighbourhoods of the form (\ref{bngen}) as well as Lemma \ref{lem-ac}, we obtain
$$
\displaystyle \frac{d\Prob^{q_1, \dots q_l}\circ T}{d{{\mathbb{P}}}^{q_1,...,q_l}}
	(X\setminus \{q_1, \dots, q_l\})=
\displaystyle \frac{d{{\mathbb{P}}}^{T{q_1},...,T{q_l}}}{d{{\mathbb{P}}}^{q_1,...,q_l}}
	(X\setminus \{q_1, \dots, q_l\}),
$$
and Proposition \ref{RN-T} follows now from
Proposition \ref{pp-cont}.
\end{proof}

\section{The relation between Palm subspaces.}

\subsection{The case of continuous measures.}
Proposition \ref{RN-T} shows that in order to establish the quasi-invariance  under the group of compactly supported diffeomorphisms for a point process, it suffices to show that its reduced Palm measures of the same order are equivalent. In this section we show that Assumption \ref{divxp} implies the relation (\ref{palm-mult-cont} between Palm subspaces; we then show that integrable kernels satisfy Assumption \ref{divxp}.
The transition from relation (\ref{palm-mult-cont} to the equivalence of  Palm measures is achieved in the next section (cf.  Corollary \ref{palm-meas-cont} and the formula (\ref{palm-meas-discr})).

As before,  we consider a closed subspace $L\subset L_2(E, \mu)$ such that  that the corresponding orthogonal projection $\Pi$ is locally trace class with a  kernel satisifying Assumptions \ref{ass-ker} , \ref{divxp}. Assume that the measure $\mu$ satisfies $\mu(\{p\})=0$ for any $p\in {\mathbb R}$.
\begin{proposition}\label{palm-subsp-cont}
For any distinct points $p_1, \dots, p_l, {q}_1,\dots, {q}_l\in U$
we have
\begin{equation}\label{palm-mult-cont}
L(p_1, \dots, p_l)=\displaystyle \displaystyle \frac{(x-p_1)\dots (x-p_l)}{(x-{q}_1)\dots (x-{q}_l)}L({q}_1, \dots, {q}_l).
\end{equation}
\end{proposition}
{\bf Remark.} The coincidence of subspaces is understood as coincidence of subspaces in $L_2$;  the functions from the right-hand side subspace are of course not defined at the points ${q}_1, \dots, {q}_l$; they are nonetheless well-defined as elements of $L_2$  since the measure $\mu$ is continuous.
For discrete measures the formulation will be modified.

\begin{proof}  In the continuous case, Assumption \ref{divxp} implies, for any $q\in U$, the inclusion
$
\displaystyle \frac{L(q)}{x-q}\subset L,
$
whence
$
\displaystyle \frac{x-p_1}{x-{q}_1}L({q}_1)\subset L({q}_1)+\frac{L(q_1)}{x-q_1}\subset L,
$
Any  function $\varphi\in \displaystyle \frac{x-p_1}{x-{q}_1}L({q}_1)$ satisfies $\varphi(p_1)=0$, whence
$
\displaystyle \frac{x-p_1}{x-{q}_1}L({q}_1)\subset L({p}_1).
$
Interchanging the points $p_1$ and ${q}_1$, we obtain the converse inclusion (using again continuity of the measure $\mu$), and the proposition is proved for $l=1$.
Proposition \ref{p1-pl} implies that if $\Pi$ satisfies Assumption \ref{divxp}, then so does $\Pi^p$ for any $p\in U$, and the proof is completed  by induction on $l$.
\end{proof}
\subsection{The case of discrete measures.}
We now let $E$ be a countable discrete subset of ${\mathbb R}$ endowed with the counting measure $\mu$.
\begin{proposition}\label{palm-subsp-discr}
Let the kernel $\Pi$ satisfy Assumption \ref{divxp}. Let $p_1, \dots, p_l\in E$ be distinct, and let $\pi$ be a permutation of $\{1, \dots, l\}$. Then we have
\begin{multline}\label{palmeqdisc}
L( p_{\pi(1)}, \dots, p_{\pi(m)}, \breve{p}_{\pi(m+1)}, \dots, \breve{p}_{\pi(l)})=\\=
\chi_{_{E\setminus \{p_1, \dots, p_l\}}}(x)\displaystyle \frac{(x-p_{\pi(1)})\dots (x-p_{\pi(m)})}{(x-{ p}_1)\dots (x-{p}_m)}
L(p_1, \dots, p_m, \breve{p}_{m+1}, \dots, \breve{p}_l).
\end{multline}
\end{proposition}

\begin{proof} As in the continuous case, we proceed by induction and start with the case $l=2, m=1$: we need to show, for any
distinct $p,q\in E$, the equality
\begin{equation}\label{palm-disc-1}
L(p, {\breve q})= \displaystyle \chi_{E\setminus \{p,q\}}\frac{x-p}{x-q}L(q, {\breve p}).
\end{equation}
Now, by Proposition \ref{lprim}, we have
$$
\chi_{E\setminus q}\displaystyle \frac {L(q)}{x-q}\subset \chi_{E\setminus q}L.
$$
Since
$$ \displaystyle \frac{x-p}{x-q} = 1 + \displaystyle \frac{q-p}{x-q},$$ we also  have
$$
	\displaystyle \frac{x-p}{x-q}\chi_{E\setminus q}L(q) \subset \chi_{E\setminus q} L = L(\breve{q})
$$
Now, multiplying any function by $\chi_{E\setminus q}\displaystyle \frac{x-p}{x-q}$ yields a function that assumes value 0 at the point $p$; we thus conclude
 \begin{equation}\label{pq-base}
\chi_{E\setminus \{p,q\}}\displaystyle \frac{x-p}{x-q}L({\breve p},q)  \subset L(p, \breve{q}).
\end{equation}
Interchanging the variables $p,q$, we obtain the inverse inclusion, and (\ref{palm-disc-1}) is proved.
If $\Pi$ satisfies Assumption \ref{divxp}, then, by Proposition \ref{p1-pl}, so does $\Pi^p$, and, by definition, so does $\Pi^{\breve p}$. Induction completes the proof. \end{proof}

\subsection{ Proof of Proposition \ref{intxpred}.}
\subsubsection{The subspace $L^{\prime}$.}
\begin{proposition}
If $\varphi\in L_2({\mathbb R}, \mu)$ is such that $x\varphi\in L_2({\mathbb R}, \mu)$, then
the integrals
$$
\displaystyle \int\limits_{\mathbb R} \varphi(x)A(x) d\mu(x), \displaystyle \int\limits_{\mathbb R} \varphi(x)B(x) d\mu(x)
$$
are well-defined.
\end{proposition}
\begin{proof}
Since for any $p\in U$ we have
$$
\langle v_p(x), (x-p)\varphi(x)\rangle =A(p) \displaystyle \int\limits_{\mathbb R} \varphi(x)B(x) d\mu(x)- B(p)\displaystyle \int\limits_{\mathbb R} \varphi(x)A(x) d\mu(x),
$$
keeping in mind that linear independence of $A$ and $B$ implies the existence of points $p_1, p_2\in U$ such that
$A(p_1)B(p_2)-A(p_2)B(p_1)\neq 0$, we obtain the desired integrability.
\end{proof}
Let
\begin{equation}
L^{\prime}=\{\psi\in L: x\psi\in  L_2({\mathbb R}, \mu),  \int\limits_{\mathbb R}\psi(x)A(x) d\mu(x)=  \int\limits_{\mathbb R}\psi(x)B(x) d\mu(x)=0\}.
\end{equation}

\begin{proposition}\label{lprim}
Let $p\in U$ and $\varphi \in L$ satisfy $\varphi(p)=0$. Then there exists $\psi \in L^{\prime}$ such that
\begin{equation}
\label{phipsi}
\varphi(x) = (x-p)\psi(x).
\end{equation}
\end{proposition}

\begin{proof} It suffices to consider the case $p=0$, $A(0) = 0, B(0) \neq 0$: since the functions $A$, $B$ are linearly independent, the general case is reduced to this particular one by
a translation of ${\mathbb R}$ and a linear unimodular change of variable (\ref{chg-var-ab}).
Let $\psi^{\prime}$ be such that (\ref{phipsi}) holds (in the continuous case, such a function $\psi^{\prime}$ is unique: we simply set $\psi^{\prime}=\varphi(x)/(x-p)$ and note that $\psi^{\prime}$ is square-integrable due to Assumption \ref{ass-ker}; in the discrete case, however, there are many such functions, differing by their value at $p=0$).
Using the integrable form of our kernel and the relation $\varphi\in L$, write
\begin{multline}
\label{comm1}
(x\Pi-\Pi x)\psi^{\prime}=x\Pi\psi^{\prime}(x) - \varphi(x) =\\= A(x) \displaystyle \int\limits_{{\mathbb R}} B(y)\psi^{\prime}(y)d\mu(y) -
B(x) \displaystyle \int \limits_{{\mathbb R}} A(y)\psi^{\prime}(y)d\mu(y).
\end{multline}

Since $\varphi \in L$, $\varphi(0)=0$, $A(0)=0$, $B(0)\neq 0$, substituting $x=0$ into (\ref{comm1}) we obtain
\begin{equation}\label{apsinul}
\displaystyle \int_{\mathbb{R}} A(y)\psi^{\prime}(y)d\mu(y) = 0.
\end{equation}
Recall that by definition we have
$$
v_0(x)=\displaystyle \frac{A(x)}{x}\in L.
$$
Dividing (\ref{comm1}) by $x$ and keeping  (\ref{apsinul}) in mind, we obtain that there exists $\alpha\in {\mathbb C}$ such that
\begin{equation}
\label{comm2}
\Pi\psi^{\prime}(x)-\psi(x)- \alpha \delta_0=v_0(x) \displaystyle \int\limits_{{\mathbb R}} B(y)\psi^{\prime}(y)d\mu(y).
\end{equation}
The extra term $\alpha \delta_0$ is only necessary in the case when $\mu({0})>0$.
It follows that we have $\psi=\psi^{\prime}+\alpha\delta_0\in L$ and, consequently,  applying the commutator $x\Pi-\Pi x$ to the function $\psi$,
that we also have
$$
\displaystyle \int_{\mathbb{R}} A(y)\psi(y)d\mu(y) =\displaystyle \int\limits_{{\mathbb R}} B(y)\psi(y)d\mu(y)=0.
$$
\end{proof}

\subsection{Integrability of the Palm kernel.}

\begin{lemma}\label{palm-int}
Let $q \in U$ be such that $\Pi(q, q) \neq 0$.
Then the kernel of the operator $\Pi^q$ has the integrable form
\begin{equation}
\Pi^q(x, y) = \displaystyle \frac{A^q(x)B^q(y) - A^q(y)B^q(x)}{x-y}
\end{equation}
where
\begin{multline}
A^{q}(x) = \displaystyle \frac{A(x)B(q)-A(q)B(x)}{\sqrt{(A(q))^2+(B(q))^2}};\\
 B^{q}(x) = \displaystyle \frac{A(x)A(q)+B(x)B(q)}{\sqrt{(A(q))^2+(B(q))^2}} - \displaystyle \frac{\sqrt{(A(q))^2+(B(q))^2}(A(x)B(q)-A(q)B(x))}{\Pi(q,q)(x-q)}.
\end{multline}
\end{lemma}

\begin{proof} We first consider the case $A(q)=0, B(q)\neq 0,$ in which
$$\Pi(x,q)=\displaystyle \frac{A(x)B(q)}{x-q}
$$
and
$$
\Pi^q(x, y)=\Pi(x,y)-\displaystyle \frac{B(q)^2A(x)A(y)}{\Pi(q,q)(x-q)(y-q)}=\displaystyle \frac{A^q(x)B^q(y) - A^q(y)B^q(x)}{x-y}
$$
with
$$
A^q(x)=A(x), B^q(x)=B(x)-\displaystyle \frac{B(q)^2A(x)}{\Pi(q,q)(x-q)},
$$
as desired. Recalling that in an integrable representation of a kernel, the functions $A$ and $B$ are defined up to a unimodular change of variables (\ref{chg-var-ab}), we reduce the general case  to the particular one by a rotation
$$
A(x)\to \displaystyle \frac{A(x)B(q)-A(q)B(x)}{\sqrt{(A(q))^2+(B(q))^2}};\ 
 B(x) \to  \displaystyle \frac{A(x)A(q)+B(x)B(q)}{\sqrt{(A(q))^2+(B(q))^2}}.
$$
\end{proof}

In the discrete case, we have a dual to Proposition \ref{palm-int}.
\begin{proposition}\label{hole-int}
Let $\Pi$ be a projection operator with an integrable kernel.
Let $q \in D$ be such that $\Pi(q, q) \neq 1$.
Then the kernel of the operator $\Pi^{\breve q}$ has  integrable form
\begin{equation}
\Pi^{\breve q}(x, y) = \displaystyle \frac{A^{\breve q} (x)B^{\breve q}(y) - A^{\breve q}(y)B^{\breve q}(x)}{x-y}
\end{equation}
where $A^{\breve q}(q)=B^{\breve q}(q)=0$ and for $x\neq q, y\neq q$ we have
\begin{multline}\label{abqbreve}
A^{\breve q}(x) = \displaystyle \frac{A(x)B(q)-A(q)B(x)}{\sqrt{(A(q))^2+(B(q))^2}};\\
 B^{\breve q}(x) = \displaystyle \frac{A(x)A(q)+B(x)B(q)}{\sqrt{(A(q))^2+(B(q))^2}} +\displaystyle \frac{\sqrt{(A(q))^2+(B(q))^2}(A(x)B(q)-A(q)B(x))}{(1-\Pi(q,q))(x-q)}.
\end{multline}
\end{proposition}

\begin{proof} Direct substitution of (\ref{abqbreve}) into (\ref{pqbreve}).
\end{proof}

\section{Multiplicative functionals and regularization}

\subsection{Outline of the section}
The multiplicative functional corresponding to the function $(\frac{x-p}{x-q})^2$
might diverge (as happens, for instance, for the Airy kernel), and in this section we describe a general scheme of regularization of multiplicative functionals. Our starting point is the formula (\ref{pig-psig}), and we show that, under certain additional assumtpions, the regularized multiplicative functional ${\overline \Psi}_g$ can still be defined  if the multiplicative functional itself diverges, much in the same way in which the Hilbert-Carleman regularization of the  Fredholm determinant is defined.  Additional technical difficulties arise because we must consider functions bounded away neither from zero nor from infinity.
In the following 9 subsections  we  go back to a general phase space $(E, \mu)$
and a general locally trace class operator $\Pi$ of orthogonal projection onto a closed subspace $L\subset L_2(E, \mu)$. We use neither integrability  nor Assumption \ref{divxp}. The main results of the section are Proposition \ref{multreg-main}
giving the existence of  regularized multiplicative functionals and Proposition \ref{multreg-main-cont} establishing their continuous dependence on the function.
Our first step is the regularization of additive functionals, for which we rely on the slow growth of variance for determinantal point processes (cf. e.g. \cite{buf-gafa}).
\subsection{Regularization of additive functionals}
Let $f:E \to {\mathbb C}$ be a Borel function.
We set $S_f$ to be the corresponding additive functional, and, if
$S_f\in L_1(\Conf(E), \Prob_{\Pi})$, then we set
\begin{equation}\label{def-overlinesf}
{\overline S}_f=S_f-\ee S_f.
\end{equation}
The random variable ${\overline S}_f$ will be called the {\it normalized} additive functional corresponding to $f$.
We shall now see that the normalized additive functional can be defined even when the additive functional itself is not well-defined.
Set
$$
\Var(\Pi, f)=\displaystyle \frac12 \displaystyle \int\limits_E\displaystyle \int\limits_E |f(x)-f(y)|^2 |\Pi(x,y)|^2d\mu(x)d\mu(y).
$$
Note that the value $\Var(\Pi, f)$ does not change if the function $f$ is changed by an additive constant.
If $S_f\in L_2(\Conf(E), \Prob_{\Pi})$, then $\Var(\Pi,f)<+\infty$ and
\begin{equation}
\Var(S_f)=\ee |{\overline S}_f|^2=\Var(\Pi,f).
\end{equation}

Note also the clear inequality
\begin{equation}\label{varin}
\Var(\Pi,f)\leq 2\displaystyle \int\limits_E |f(x)|^2 \Pi(x,x)d\mu(x)
\end{equation}
which is obtained by summing the inequality $|f(x)-f(y)|^2\leq 2(|f(x)|^2+|f(y)|^2)$
over all $x,y$ and using the Pythagoras theorem (reproducing property of $\Pi$ on the diagonal):
$$
\Pi(x,x)=\displaystyle \int\limits_{ E} |\Pi(x,y)|^2 d\mu(y).
$$

The integral defining the variance of an additive functional may converge even when the integral defining its expectation does not: for instance, if $f(x)=(|x|+1)^{-1}$ and $\Pi$ is the discrete sine-kernel. The normalized additive functional can nonetheless by continuity be defined in $L_2$ even when the additive functional itself diverges almost surely.

Introduce the Hilbert space ${\mathcal V}(\Pi)$ in the following way: the elements of  ${\mathcal V}(\Pi)$ are functions $f$ on $E$ satisfying $\Var(\Pi, f)<+\infty$; functions that differ by a constant are identified, but, slightly abusing terminology we still refer to elements of ${\mathcal V}(\Pi)$ as functions.
The square of the norm of an element $f\in {\mathcal V}(\Pi)$ is precisely $\Var(\Pi,f)$.
By definition, bounded functions that are identically zero in the complement of a bounded set form a dense subset of
${\mathcal V}(\Pi)$. The correspondence $f\to {\overline S}_{f}$ is thus an isometric embedding of a dense subset of ${\mathcal V}(\Pi)$ into $L_2(\Conf(E), \Prob_{\Pi})$; it  therefore admits a unique  isometric extension onto the whole space ${\mathcal V}(\Pi)$, and we  obtain the following
\begin{proposition}
There exists a unique linear isometric embedding
$${\overline S}: {\mathcal V}(\Pi)\to  L_2(\Conf(E), \Prob_{\Pi}), \
{\overline S}: f\to {\overline S}_f
$$
such that
\begin{enumerate}
\item $\ee {\overline S}_f=0$ for all $f\in {\mathcal V}(\Pi)$;
\item if $S_f\in L_1(\Conf(E), \Prob_{\Pi})$, then ${\overline S}_f$ is given by (\ref{def-overlinesf}).
\end{enumerate}

\end{proposition}

\subsection{Regularization of multiplicative functionals.}

Given a function $g$ such that $\Var(\Pi,\log g)<+\infty$, set
$$
{\widetilde \Psi}_g=\exp({\overline S}_{\log g}).
$$
By definition, we have
\begin{equation}\label{multiplicat}
{\widetilde \Psi}_{g_1g_2}={\widetilde \Psi}_{g_1}{\widetilde \Psi}_{g_2}.
\end{equation}
Since $\ee {\overline S}_{\log g}=0$, by Jensen's inequality, for any positive function $g$ we  have
$$
\ee {\widetilde \Psi}_g\geq 1.
$$
The expectation $\ee {\widetilde \Psi}_g$ may however be infinite, and our next aim is to give conditions for its finiteness.

It will be convenient for us to allow zero values for the function $g$: let therefore $g$ be nonnegative,
set $E_0=\{x\in E: g(x)=0\}$,  assume that the subset
$\Conf(E; E\setminus E_0)$ of those configurations all whose particles lie in $E\setminus E_0$ has positive probability, consider the restriction of our measure $\Prob$ onto the subspace
$\Conf(E; E\setminus E_0)$, introduce the corresponding functional ${\widetilde \Psi}_g$ and extend it to the whole of $E$ by setting
${\widetilde \Psi}_g(X)=0$ for all configurations containing a particle at $E_0$.
Assume that $\mathrm{tr} \chi_{E_0}\Pi\chi_{E_0}<+\infty$.
Then we have  $\Prob_{\Pi}(\Conf(E; E\setminus E_0))=\det (1-\chi_{E_0}\Pi\chi_{E_0})$. The operator $\chi_{E_0}\Pi\chi_{E_0}$ is positive, contractive and has finite trace, so the operator $1-\chi_{E_0}\Pi\chi_{E_0}$ is invertible once it is injective.
Thus, if a function
$\varphi\in L$ satisfying $\varphi(x)=0$ for all  $x\in E\setminus E_0$ must be the zero function, then
 $\Prob_{\Pi}(\Conf(E; E\setminus E_0))>0$.

If ${\widetilde \Psi}_g\in L_1(\Conf(E), \Prob_{\Pi})$, then, as before,  we write
$$
{\overline \Psi}_g=\displaystyle \frac{{\widetilde \Psi}_g}{\ee{\widetilde \Psi}_g}.
$$
The main result of this section is
\begin{proposition}\label{multreg-main}

Let $g$ be a nonnegative function
 satisfying the following:
\begin{enumerate}
\item there exist $\varepsilon>0$, $M>\varepsilon$ such that  the set
$E^{\varepsilon, M}=\{x\in E: g(x)<\varepsilon \ \mathrm{or} \ g(x)>M \}$ is bounded and
$$
|| \chi_{E^{\varepsilon, M}}\Pi ||<1;
$$
\item $$
\int\limits_{E^{\varepsilon, M}} |g(x)|\Pi(x,x)d\mu(x)+ \int\limits_{E\setminus E^{\varepsilon, M}} |g(x)-1|^2\Pi(x,x) d\mu(x)<+\infty.
$$
\end{enumerate}
Then ${\widetilde \Psi}_g\in L_1(\Conf(E), \Prob_{\Pi})$, the subspace $\sqrt{g}L$ is closed, the corresponding operator of orthogonal projection $\Pi^g$ is locally of trace class, and we have
\begin{equation}\label{probpig1}
\Prob_{\Pi^g}={\overline \Psi}_g\Prob_{\Pi}.
\end{equation}
\end{proposition}

We also need to establish continuity of ${\overline \Psi_g}$ as a function of $g$.
This continuity is established in a specially constructed function space  whose definition is a bit involved.
We fix positive  numbers $\alpha>0$, $\varepsilon>0$, $M>\varepsilon$ and two bounded Borel subsets $B^1$, $B^2$  of $E$ satisfying
$$
|| \chi_{B^1 \cup B^2}\Pi ||<1.
$$

We now let $\mathscr G$ be the set of nonnegative measurable functions $g$  on $E$ satisfying
\begin{enumerate}
\item $\{x: g(x)<\varepsilon\}\subset B^1$;
\item $\{x: g(x)>M\}\subset B^2$;
\item  $\int\limits_{B^2} |g(x)|^{1+\alpha}\Pi(x,x)d\mu(x)+ \int\limits_{E\setminus B^2} |g(x)-1|^2\Pi(x,x) d\mu(x)<+\infty$.
\end{enumerate}
We metrize the set $\mathscr G$ by setting
$$
d_{{\mathscr{G}}}(g_1, g_2)=\int\limits_{B^2} |g_1(x)-g_2(x)|^{1+\alpha}\Pi(x,x)d\mu(x)+ \int\limits_{E\setminus B^2} |g_1(x)-g_2(x)|^2\Pi(x,x) d\mu(x).
$$
The distance $d_{{\mathscr{G}}}$ turns $\mathscr G$ into a complete separable metric space.

\begin{proposition}\label{multreg-main-cont}
For any $\alpha^{\prime}: 0<\alpha^{\prime}<\alpha$, the correspondences $g\to \widetilde \Psi_g$,  $g\to{\overline \Psi}_g$
induce continuous mappings from $\mathscr G$ to $L_{1+\alpha^{\prime}}(\Conf(E), \Prob_{\Pi})$,
\end{proposition}
{\bf Remark.} Along similar lines, one also proves that the correspondence
$g\to \Pi^g$ induces a continuous mapping from the space $\mathscr G$ to
$\mathscr{I}_{1, \mathrm{loc}}(E, \mu)$.

Propositions \ref{regmult-cont}, \ref{regmult-discr} directly follow from  Propositions \ref{multreg-main}, \ref{multreg-main-cont}.

The proof of the main propositions \ref{multreg-main} and \ref{multreg-main-cont} occupies the bulk of this section. We separately consider the ``main part'' (where the function $g$ is bounded away from $0$ and $\infty$),
the neighbourhood of $0$ and the neighbourhood of $\infty$. The proof is concluded by decomposing a general function $g$ as a product of these three types of functions.

First, we reduce ourselves to the case of {\it positive} $g$. Let $E_0=\{x: g(x)=0\}$.
 Our assumptions imply $\Prob_{\Pi}(\Conf(E; E\setminus E_0))>0$ and so,
 restricting ourselves, if necessary,  to the subset $\Conf(E; E\setminus E_0)$, we can assume that the function $g$ is positive.

{\bf Remark.} Our aim is to apply Proposition \ref{multreg-main} to functions of the form $g(x)= \frac{(x-p_1)\dots (x-p_l)}{(x-{q}_1)\dots (x-{q}_l)}$. In the continuous case, such functions are almost surely nonzero. In the discrete case, however,
the finite zero set of our function  has positive measure: whence the need for the extra subset $E_0$ in our proposition,
and the need of the assumption that the subspace $L$ not admit finitely-supported functions in the formulation of
Theorem \ref{main-discr}.

\subsection{An estimate of diagonal values of the kernel $\Pi^g$.}

Iterating multiplicative functionals, we need to estimate traces of the resulting kernels.  The following simple proposition shows that  diagonal values of the kernel of $\Pi^g$ can be estimated from above by the diagonal values of the kernel $\Pi$.

\begin{proposition}\label{kernelest}
Let the kernel $\Pi$ satisfy Assumption \ref{ass-ker}, and let $g$ be a non-negative bounded function on $E$ such that the operator $1+(g-1)\Pi$ is invertible. Then for any $q\in U$ we have
\begin{equation}
\label{pigdiag}
\Pi^{g}(q,q)\leq  {g(q)} ||(1+(g-1)\Pi)^{-1}||\Pi(q,q).
\end{equation}
\end{proposition}
\begin{proof} As before, we let $\langle, \rangle$ be the standard inner product in $L_2(E, \mu)$
and we write $v_q(x)=\Pi(x,q)$ so that $\Pi(q,q)=\langle v_q, v_q\rangle$. By definition then
$$
\Pi^{g}(q,q)=g(q)  \langle \Pi(1+(g-1)\Pi)^{-1}v_{q}, v_q\rangle\leq  {g(q)} ||(1+(g-1)\Pi)^{-1}||\langle v_q, v_q\rangle.
$$
\end{proof}

\subsection{The case of functions $g$ bounded away from $0$ and $\infty$.}

Let $\A_2(\Pi)$ be the set of positive Borel functions $g$ on $E$ satisfying
\begin{enumerate}
\item
$0< \inf\limits_E g\leq \sup\limits_E g<\infty$;
\item
$$
\displaystyle \int\limits_E |g(x)-1|^2\Pi(x,x)d\mu(x)<+\infty.
$$
\end{enumerate}

By definition, the set $\A_2(\Pi)$ is a semigroup under multiplication.

Endow the set $\A_2(\Pi)$  with a metric by setting the distance between two functions $g_1$ and $g_2$ to be
$$
\sqrt{\displaystyle \int\limits_E |g_1(x)-g_2(x)|^2\Pi(x,x)d\mu(x)}.
$$

Using  the second condition in the definition of $\A_2(\Pi)$ and the estimate (\ref{varin}), for any  $g\in\A_2(\Pi)$ we have
$$
\Var(\Pi, g-1)<+\infty.
$$

Since on any interval of the positive half-line, bounded away from zero and infinity, the quantity
$
|\log t-t+1|/t^2
$
is bounded both above and below,
for any function $g\in \A_2(\Pi)$, we also have
$$
\Var(\Pi, \log g)<+\infty.
$$
In particular, for any function $g\in \A_2(\Pi)$ the functional ${\widetilde \Psi}_g$ is well-defined.
The following proposition, the main result of this subsection, establishes  its integrability.
\begin{proposition}\label{cont-bnd}
For any  $p\geq 1$, and any function $g\in \A_2(\Pi)$ we have ${\widetilde \Psi}_g\in L_p(\Conf(E), \Prob_{\Pi})$.
The correspondences
$$
g\to {\widetilde \Psi}_g, g\to {\overline \Psi}_g
$$
are continuous mappings from $\A_2(\Pi)$ to $L_p(\Conf(E), \Prob_{\Pi})$.
\end{proposition}
\begin{proof} Let \begin{equation}\label{meps}
\A^{\varepsilon, M}_2(\Pi)=\{g\in \A_2(\Pi): \varepsilon \leq \inf\limits_E g\leq \sup\limits_E g\leq M\}.
\end{equation}
\begin{proposition}\label{l2-est}
For any $\varepsilon>0, M>0$ there exists a constant
$C_{\varepsilon,M}>0$ such that any $g\in \A^{\varepsilon, M}_2(\Pi)$ satisfies
\begin{equation}\label{l2}
 \log\ee  |{\widetilde \Psi}_{g}|^2\leq C_{\varepsilon, M} \displaystyle \int\limits_E |g(x)-1|^2\Pi(x,x)d\mu(x);
\end{equation}
\end{proposition}
\begin{proof} It suffices to prove the estimate
\begin{equation}\label{l1}
\log\ee  {\widetilde \Psi}_{g}\leq C_{\varepsilon, M} \displaystyle \int\limits_E |g(x)-1|^2\Pi(x,x)d\mu(x),
\end{equation}
and (\ref{l2}) follows by multiplicativity (\ref{multiplicat}).
It suffices to establish (\ref{l1}) in the case when the set $\{x\in E: g(x)\neq 1\}$ is bounded, as the general case
follows by  Fatou's lemma.
Now there exists a constant $C_2>0$ such that
\begin{multline}
\log \ee \Psi_g\leq \tr (\sqrt{g-1}\Pi\sqrt{g-1})+C_2||\sqrt{g-1}\Pi\sqrt{g-1}||_2^2= \\
=\displaystyle \int\limits_E (g(x)-1) \Pi(x,x) d\mu(x)+C_2\displaystyle \int\limits_E |g(x)-1|^2\Pi(x,x)d\mu(x).
\end{multline}
We have assumed boundedness of the set $\{x\in E: g(x)\neq 1\}$ in order that the integral
$
\displaystyle \int\limits_E (g(x)-1) \Pi(x,x) d\mu(x)
$
be well-defined; it will, however, disappear from the final result.
From (\ref{meps}), again using the fact that the quantity
$
|\log t-t+1|/t^2
$
is bounded both above and below by constants only depending on $\varepsilon$ and $M$, we obtain
\begin{multline}
\left|\displaystyle \int\limits_E (g(x)-1) \Pi(x,x) d\mu(x) -\displaystyle \int \limits_E \log g(x) \Pi(x,x) d\mu(x)\right|\leq \\ \leq C_{\varepsilon, M}\displaystyle \int\limits_E |g(x)-1|^2\Pi(x,x)d\mu(x),
\end{multline}
whence finally
$$
\log \ee {\widetilde \Psi}_g=\log \ee \Psi_g-\ee S_{\log g}\leq C^{\prime}_{ \varepsilon, M}\displaystyle \int\limits_E |g(x)-1|^2\Pi(x,x)d\mu(x)
$$
\end{proof}

\begin{proposition}\label{l1-est}
For any $\varepsilon>0, M>0$ there exists a constant
$C_{\varepsilon,M}>0$ such that for $g_1, g_2\in \A_2^{\varepsilon, M}(\Pi)$ we have
$$
\ee | {\widetilde \Psi}_{g_1}- {\widetilde \Psi}_{g_2}|\leq \ee|{\widetilde \Psi}_{g_1}|^2\left(\exp\left(C_{\varepsilon,M}\ \displaystyle \int\limits_E  |g_1(x)-g_2(x)|^2\Pi(x,x)d\mu(x)\right)-1\right).
$$
For any $p\geq 1$, the correspondence $g\to {\widetilde \Psi}_g$ induces a continuous mapping
from $\A_2^{\varepsilon, M}(\Pi)$ to $L_p(\Conf(E), \Prob_{\Pi})$.
\end{proposition}

\begin{proof} Since $\ee  {\widetilde \Psi}_{g}\geq 1$, we have
$$
 \ee | {\widetilde \Psi}_{g}-1|^2\leq \ee {\widetilde \Psi}_{g^2}-1.
$$
From
the estimate (\ref{l2}) we have
\begin{equation}\label{contat1}
\ee | {\widetilde \Psi}_{g}-1|^2\leq \exp\left(C\displaystyle \int\limits_E  |g(x)-1|^2\Pi(x,x)d\mu(x)\right)-1.
\end{equation}
Applying (\ref{contat1}) to $g=g_1/g_2$, recalling the boundedness of both $g_1$ and $g_2$ and using multiplicativity, we obtain the proposition. Since, for any $p\geq 1$,  we have  ${\widetilde \Psi}_{g^p}=\left({\widetilde \Psi}_g\right)^p$,
and $\left(\A_2^{\varepsilon, M}(\Pi)\right)^p\subset \A_2^{\varepsilon^p, M^p}(\Pi)$, the desired continuity in $L_p$ follows as well.
\end{proof}

Proposition \ref{l1-est} implies Proposition \ref{cont-bnd}. \end{proof}

\subsection{The  Radon-Nikodym derivative}
\begin{corollary}\label{ga2}
Let $g\in \A_2(\Pi)$ be such that the operator $1+(g-1)\Pi$ is invertible.
Then the operator $\Pi^g$ is locally of trace class, and
 we have
\begin{equation}\label{probpig2}
\Prob_{\Pi^g}={\overline \Psi}_g\Prob_{\Pi}.
\end{equation}
\end{corollary}
\begin{proof}
Let $E^{(n)}$ be a sequence of bounded sets exhausting $E$, set $g_n=1+(g-1)\chi_{E^{(n)}}$, and note that
$||(g_n-g)\Pi||\to 0$ as $n\to\infty$. For any $n$, we have
$$
\Prob_{\Pi^{g_n}}={\overline \Psi}_{g_n}\Prob_{\Pi}.
$$
The operators
$
\Pi^g=\sqrt{g}\Pi(1+(g-1)\Pi)^{-1}\sqrt{g},
$
$
\Pi^{g_n}=\sqrt{g_n}\Pi(1+(g_n-1)\Pi)^{-1}\sqrt{g_n}
$ are locally of trace class since so is $\Pi$.
Proposition \ref{pign-conv} implies that
$
\Pi^{g_n}\to \Pi^g
$
in $\scrI_{1, \text{loc}}(E,\mu)$, and
Proposition \ref{loctracetoweak}
implies that, as $n\to\infty$,  the sequence of
measures $\Prob_{\Pi^{g_n}}$ weakly converges to $\Prob_{\Pi^g}$
in the space of probability measures on $\Conf(E)$.
Proposition \ref{cont-bnd} implies ${\overline \Psi}_{g_n}\to {\overline \Psi}_g$ in $L_1(\Conf(E), \Prob_{\Pi})$, whence
${\overline \Psi}_{g_n}\Prob_{\Pi}\to {\overline \Psi}_g\Prob_{\Pi}$
weakly in the space of probability measures on $\Conf(E)$, implying
\eqref{probpig2}.
\end{proof}

\subsection{Multiplicative functionals corresponding to a function  $g\geq 1$.}
\subsubsection{The case of bounded $g$.}

Proposition \ref{mult-prelim} takes a simpler form when our bounded  function $g$ satisfies $g\geq 1$. First, in this case the subspace $\sqrt{g}L$ is automatically closed. Second, if
$\sqrt{g-1}\Pi\sqrt{g-1}$ belongs to the trace class, then the operator $1+(g-1)\Pi$ is automatically
invertible. To verify this, observe first that in this case the operator $\sqrt{g-1}\Pi$ is Hilbert- Schmidt,
consequently, the operator $(g-1)\Pi$ is also Hilbert-Schmidt and, a fortiori, compact. To check the invertibility of the operator $1+(g-1)\Pi$, it thus suffices to check its injectivity, that is, to
to prove that a function $\varphi$ satisfying
\begin{equation}\label{varphip}
\varphi+(g-1)\Pi\varphi=0
\end{equation}
 must be the zero function.
Set $\psi=-\sqrt{g-1}\Pi\varphi$ so that $\varphi=\sqrt{g-1}\psi$. By definition, both $\varphi$ and $\psi$ are  zero on the set $\{x\in E: g(x)=1\}$. From (\ref{varphip}) we now have
$$
\psi+\sqrt{g-1}\Pi\sqrt{g-1}\psi=0,
$$
whence
$$
\langle \psi, \psi\rangle + \langle \Pi \varphi, \varphi\rangle=0,
$$
whence finally $\varphi=\psi=0$.

We can now reformulate Proposition \ref{mult-prelim} in the following simpler form

\begin{proposition}\label{g1bd}
Let $g$ be a bounded measurable function on $E$ satisfying $g\geq 1$ and such that the operator
$\sqrt{g-1}\Pi\sqrt{g-1}$ belongs to the trace class. Then all  the conclusions of Proposition \ref{mult-prelim} hold for the function $g$.
\end{proposition}

\subsubsection{The case of unbounded $g$.}

The function $(x-p)/(x-q)$ is unbounded on ${\mathbb R}$, and we prepare, for future use, a proposition
on multiplicative functionals corresponding to unbounded functions.
As before, we start with a locally trace-class operator $\Pi$ of orthogonal projection onto a subspace $L$.
We consider a function $g\geq 1$  such that the space
$\sqrt{g}L$ is a subspace of $L_2(E, \mu)$; since $g\geq 1$, it is automatically a closed subspace, and we let $\Pi^g$ be the operator of orthogonal projection onto $\sqrt{g}L$.

\begin{proposition}\label{g1unbd}

Let $\alpha>0$ and let $g$ be a Borel function on $E$ satisfying $g>1+\alpha$. Assume that
\begin{enumerate}
\item we have $\sqrt{g}L\subset L_2(E, \mu)$;
\item the operator $\sqrt{g-1}\Pi$ is Hilbert-Schmidt.
\end{enumerate}
Then $\Psi_g\in L_1(\Conf(E), \Prob_{\Pi})$, and
we have
\begin{equation}\label{pig-psig-3}
\displaystyle \frac {{\Psi}_g\Prob_{\Pi}}{\int\limits_{\Conf(E)} {\Psi}_gd\Prob_{\Pi} }=\Prob_{\Pi^g}.
\end{equation}
\end{proposition}
\begin{proof}
For $R>0$ set $g^R(x)=g(x)$ if $g(x)<R$ and $g^R(x)=1$ otherwise. As in Section 2, we use the notation
$\Pi_g=\sqrt{g-1}\Pi\sqrt{g-1}$.

The Hilbert-Schmidt norm of the operator $\sqrt{g-1}\Pi$ is given by the formula
\begin{equation}\label{hsg}
||\sqrt{g-1}\Pi||_2=\displaystyle \int\limits_E (g(x)-1)\Pi(x,x)d\mu(x).
\end{equation}
For $R_1>R_2>0$, we have
\begin{equation}\label{r1r2}
||\sqrt{g^{R_1}-1}\Pi-\sqrt{g^{R_2}-1}\Pi||_2=\int\limits_E  (\sqrt{g^{R_1}(x)-1}-\sqrt{g^{R_2}(x)-1})^2\Pi(x,x)d\mu(x),
\end{equation}
and, in view of the convergence of the integral in (\ref{hsg}),
 the right-hand side  of (\ref{r1r2}) becomes arbitrarily small once $R_1, R_2$ are sufficiently large.
It follows that, as $R\to\infty$, the sequence of operators $\sqrt{g^R-1}\Pi$ converges in Hilbert-Schmidt norm, and, consequently, the sequence of operators $\Pi_{g^R}$ converges in the trace-class norm.
Write
$
\Pi_g=\lim\limits_{R\to\infty} \Pi_{g^R},
$
and note that trace-class convergence implies
\begin{equation}\label{detpig}
\det(1+\Pi_g)=\lim\limits_{R\to\infty} \det(1+\Pi_{g^R}).
\end{equation}
For any $X$, as $R$ grows, the sequence $\Psi_{g^R}(X)$  increases (possibly assuming infinite value starting from some $R$).
By definition of the multiplicative functional, we have pointwise convergence
\begin{equation}\label{psigrpsig}
\Psi_g=\lim\limits_{R\to\infty} \Psi_{g^R},
\end{equation}
and, by (\ref{detpig}), this convergence also takes place in $L_1(\Conf(E), \Prob_{\Pi})$.
We thus finally have
$$\int\limits_{\Conf(E)} \Psi_g d\Prob_{\Pi}=\det(1+\Pi_g).$$
Using  \eqref{pig-pig},  write
$$
\Pi^{g^R}=\sqrt{1+(g^R-1)^{-1}}\Pi_{g^R}(1+\Pi_{g^R})^{-1}\sqrt{1+(g^R-1)^{-1}}.
$$
It follows that the sequence of operators $\Pi^{g^R}$ converges, as $R\to\infty$, in the space of locally trace-class operators. Since $\Pi^{g^R}$ is the operator of orthogonal projection onto the subspace $\sqrt{g^R}L$, we obtain
that
\begin{equation}\label{pigrpig}
\Pi^g=\lim\limits_{R\to\infty} \Pi^{g^R} \  \mathrm{in}  \ \mathscr I_{1,  \mathrm{loc}}(E,\mu).
\end{equation}
Now, for any fixed $R$, by Proposition \ref{g1bd}, we have
\begin{equation}\label{pig-psig-3-prelim}
\displaystyle \frac {{\Psi}_{g^R}\Prob_{\Pi}}{\int\limits_{\Conf(E)} {\Psi}_{g^R}d\Prob_{\Pi} }=\Prob_{\Pi^{g^R}}.
\end{equation}
Proposition \ref{loctracetoweak} implies that the sequence of measures $\Prob_{\Pi^{g^R}}$ weakly converges to $\Prob_{\Pi^g}$ as $R\to\infty$,
and (\ref{pig-psig-3-prelim}), together with \eqref{psigrpsig} and \eqref{pigrpig}, implies the desired relation
 (\ref{pig-psig-3}).
\end{proof}

\subsection{Conclusion of  the proof of Proposition \ref{multreg-main}.}
Set
\begin{equation}\label{def-g1}
g_0=(g-1)\chi_{E\setminus E^{\varepsilon, M}}+1.
\end{equation}
\begin{equation}\label{def-g2}
g_1=(g-1)\chi_{\{x\in E: g(x)<\varepsilon\}}+1.
\end{equation}
\begin{equation}\label{def-g3}
g_2=(g-1)\chi_{\{x\in E: g(x)>M\}}+1.
\end{equation}
By definition, $g=g_0g_1g_2$.

By definition, the subspace $\sqrt{g_0}L$ is closed, and,
by Corollary \ref{ga2}, we have
$
\Prob_{\Pi^{g_0}}={\overline \Psi_{g_0}}\Prob_{\Pi}.
$
Proposition \ref{kernelest} implies the existence of a positive constant $C$
such that
$
\Pi^{g_0}(x,x)\leq C\Pi(x,x)
$
for $\mu$-amost all $x\in E$.
Applying Proposition \ref{gless1bd} to the function $g_1$ and the operator $\Pi^{g_0}$,  we arrive at the formula
$$
\Prob_{\Pi^{g_1g_0}}={\overline \Psi_{g_1}}\Prob_{\Pi^{g_0}}={\overline \Psi_{g_1g_0}}\Prob_{\Pi}.
$$
Again, Proposition \ref{kernelest} implies the existence of a positive constant $C$
such that $
\Pi^{g_1g_0}(x,x)\leq C\Pi(x,x)
$
for $\mu$-amost all $x\in E$, whence
$$
\displaystyle \int\limits_E |g_2(x)-1|\Pi^{g_1g_0}(x,x)d\mu(x)<+\infty,
$$
and, consequently, that the operator $\sqrt{g_2-1}\Pi^{g_1g_0}$ is Hilbert-Schmidt.

We therefore apply Proposition \ref{g1unbd} to the function $g_2$ and the operator $\Pi^{g_1g_0}$ and obtain
$$
\Prob_{\Pi^{g_2g_1g_0}}={\overline \Psi_{g_2}}\Prob_{\Pi^{g_1g_0}}.
$$
Observe that we only used regularized multiplicative functionals at the very first step of our argument.
In other words, there exist constants $C_1$, $C_2$, $C_3$ such that we have
$$
\Prob_{\Pi^{g_0}}=C_0{\widetilde \Psi}_{g_0}\Prob_{\Pi}, \
\Prob_{\Pi^{g_1g_0}}=C_1\Psi_{g_1} \Prob_{\Pi^{g_0}}, \
\Prob_{\Pi^g}=\Prob_{\Pi^{g_2g_1g_0}}=C_2\Psi_{g_2}\Prob_{\Pi^{g_1g_0}}.
$$
By definition, we have
$
{\widetilde \Psi}_g={\widetilde \Psi}_{g_2}{\widetilde \Psi}_{g_1}{\widetilde \Psi}_{g_0}
$
and, consequently, for a suitable positive constant $C^{\prime}$, also
$
{\widetilde \Psi}_g=C^{\prime}{\widetilde \Psi}_{g_0}{\Psi}_{g_1}{\Psi}_{g_2}.
$
Passing from ${\widetilde \Psi}_g$ to ${\overline \Psi}_g$ and noting that the Radon-Nikodym derivative must by definition have expectation $1$,  we finally obtain the desired relation
$
\Prob_{\Pi^g}={\overline \Psi}_g\Prob_{\Pi}.
$
\qed

\subsection{Proof of Proposition \ref{multreg-main-cont}.}
Set $\mathscr{G}^0=\{g\in \mathscr{G}: (g_0-1)\chi_{E\setminus (B^1\cup B^2)}=g_0-1\}$ and endow the set $\mathscr{G}^0$ with topology induced from $\mathscr{G}$, or, equivalently, from $L_2(E, \mu)$.
Set $\mathscr{G}^1=\{g_1\in \mathscr{G}: (g_1-1)\chi_{B^1}=g_1-1 \}$, and endow the set $\mathscr{G}^1$ with topology induced from $\mathscr{G}$, or, equivalently, from $L_2(E, \mu)$ or, equivalently again, from $L_1(E, \mu)$ (recall here that our functions are uniformly  bounded.
Finally, set $\mathscr{G}^2=\{g\in \mathscr{G}: (g_2-1)\chi_{ B^2}=g_2-1\}$
endow the set $\mathscr{G}^2$ with topology induced from $\mathscr{G}$, or, equivalently, from $L_{1+\alpha}(E, \mu)$.
The continuity  of the mapping $g_0\to {\widetilde \Psi}_{g_0}$  from $\mathscr{G}^0$ to $L_p(E, \mu)$  directly follows from Propositions \ref{cont-bnd}
, \ref{l1-est}. Next, the correspondence  $g_1\to {\widetilde \Psi}_{g_1}$ is a continuous
mapping from  $\mathscr{G}^1$ to $L_p(E, \mu)$ by bounded convergence theorem.
The correspondence  $g_2\to {\Psi}_{g_2}$ induces a continuous mapping from the space
$\mathscr{G}^2$ to $L_{1+\alpha}(E, \mu)$ by the dominated convergence theorem.
Setting $g=g_0g_1g_2$,
 using the Hoelder inequality (and keeping in mind that $p$ can be taken arbitrarily large), we obtain
the desired continuity of the correspondences $g\to {\widetilde \Psi}_g$, $g\to {\overline \Psi}_g$ considered as
mappings from  $\mathscr{G}$ to $L_{1+\alpha^{\prime}}(E, \mu)$.
\qed

\subsection{Proof of Proposition \ref{regmult-cont} and Theorem \ref{main-cont}.}
Let $l\in {\mathbb N}$ and take distinct points $p_1, \dots, p_l$, ${q}_1,\dots, {q}_l\in U$. We check that one can choose $\varepsilon>0$ small enough in such a way that  the function
$$
g(x)=\left(\displaystyle \displaystyle \frac{(x-p_1)\dots (x-p_l)}{(x-{q}_1)\dots (x-{q}_l)}\right)^2
$$
satisfies the assumptions of Proposition \ref{multreg-main-cont} with an arbitrary $\alpha<3/2$ and $B^1=\{x\in {\mathbb R}: \max\limits_{i=1, \dots, l} |x-{p}_i|<\varepsilon\}$, $B^2=\{x\in {\mathbb R}: \max\limits_{i=1, \dots, l} |x-{q}_i|<\varepsilon\}$. In this proof, $C$ stands for a constant depending only on $q_1, \dots, q_l$ and $\varepsilon$.
The kernel $\Pi$ is smooth, and if $|x-q_i|<\varepsilon$, then
$
|\Pi^{{q}_1,\dots, {q}_l}(x,x)|<C|x-{q}_i|^2,
$
whence  for any
$\alpha<3/2$ we have $$\int\limits_{B^2}|g(x)|^{1+\alpha} \Pi^{{q}_1,\dots, {q}_l}(x,x)dx<+\infty.$$

 For any $\varepsilon>0$, on the set $\{x\in {\mathbb R}: \min\limits_{i=1, \dots, l} |x-{q}_i|>\varepsilon\}$ we have $|g(x)-1|^2\leq  C(1+x^2)^{-1},
$ whence
$$
\displaystyle \int\limits_{\{x\in {\mathbb R}: \min\limits_{i=1, \dots, l} |x-{q}_i|>\varepsilon\}}
 |g(x)-1|^2\Pi(x,x)d\mu(x)<+\infty,
$$
The operator $\Pi^{ {q}_1,\dots, {q}_l}$ is a finite-rank perturbation of the operator $\Pi$, whence
\begin{equation}\label{eq32}
\displaystyle \int\limits_{\{x\in {\mathbb R}: \min\limits_{i=1, \dots, l} |x-{q}_i|>\varepsilon\}}
 |g(x)-1|^2\Pi^{ {q}_1,\dots, {q}_l}(x,x)d\mu(x)<+\infty.
\end{equation}
Proposition \ref{multreg-main-cont} now implies Proposition \ref{regmult-cont}. \qed

Proposition \ref{multreg-main} together with  Proposition \ref{palm-subsp-cont} immediately imply 
\begin{corollary} \label{palm-meas-cont} Under the assumptions of Theorem \ref{main-cont}, for any distinct points $p_1, \dots, p_l, {q}_1,\dots, {q}_l\in U$, for the corresponding reduced Palm measures are equivalent, and we have
$$
\displaystyle \displaystyle \frac{{d{\mathbb P}_{\Pi^{{p}_1,\dots, { p}_l}}}}{{d{{\mathbb P}_{\Pi^{{q}_1,\dots, {q}_l}}}}}={\overline \Psi}_{\left| \frac{(x-p_1)\dots (x-p_l)}{(x-{q}_1)\dots (x-{q}_l)}\right|^2}.
$$
\end{corollary}
Together with Proposition \ref{RN-T}, Corollary \ref{palm-meas-cont} implies Theorem \ref{main-cont}.
\qed

\subsection{Proof of Proposition \ref{regmult-discr} and Theorem \ref{main-discr}.}
Denote $q_i=\sigma(p_i)$, $i=1, \dots, l$; of course, we have $\{p_1, \dots, p_l\}=\{q_1, \dots, q_l\}$. Set
\begin{equation}\label{gdisc}
g(x)=\prod\limits_{i=1}^m \left(\displaystyle \frac{x-q_i}{x-p_i}\right)^2\chi_{E\setminus \{p_1, \dots, p_l\}}(x).
\end{equation}
The function $g$ is bounded,  the condition (\ref{nsq}) implies
$$
\sum\limits_E |g(x)-1|^2 <+\infty,
$$
and, since  no nonzero function in $L$ is supported on a finite set, we conclude applying Proposition \ref{cont-bnd} to the function $g$ and the kernel $\Pi^{p_1, \dots, p_m, \breve{p}_{m+1}, \dots, \breve{p}_l}$.

In a similar way to the continuous case,
Proposition \ref{multreg-main} and Proposition \ref{palm-subsp-discr} together imply that,
under the assumptions of Proposition \ref{regmult-discr}, we have
\begin{equation}\label{palm-meas-discr}
\Prob_{\Pi^{q_1, \dots, q_m, \breve{q}_{m+1}, \dots, \breve{q}_{l}}}={\overline \Psi}(p_1, \dots, p_l, m, \sigma)
\Prob_{\Pi^{p_1, \dots, p_m, \breve{p}_{m+1}, \dots, \breve{p}_l}}.
\end{equation}
The relation (\ref{palm-meas-discr}) together with Proposition \ref{RN-T} implies Theorem \ref{main-discr}.\qed


\begin{thebibliography}{9}

\bibitem{BLPS}
I. Benjamini, R. Lyons, Y. Peres, and O. Schramm,
Uniform spanning forests, Ann. Probab.  29, 1 (2001), 1-65.




\bibitem{BOO}A. Borodin, A. Okounkov; G. Olshanski, Asymptotics of Plancherel measures for symmetric groups,
J. Amer. Math. Soc. 13 (2000), 481--515.

\bibitem{BO-gamma} A. Borodin and G. Olshanski, Random partitions and the Gamma kernel. Advances in Math. 194 (2005), 141–202; arXiv: math-ph/0305043.



\bibitem{BorRains}A.M. Borodin, E.M. Rains, Eynard-Mehta theorem, Schur process, and their pfaffian analogs.
J. Stat. Phys. 121 (2005), 291--317.

\bibitem{Buf-CIRM} A. Bufetov, Infinite determinantal measures, Electronic Research Announcements in the Mathematical Sciences, 20 (2013), pp. 8 -- 20.


\bibitem{Buf-umn}A.I. Bufetov, Multiplicative functionals of determinantal processes,
Uspekhi Mat. Nauk 67 (2012), no. 1 (403), 177--178;
translation in Russian Math. Surveys 67 (2012), no. 1, 181--182.

\bibitem{buf-gafa} A. I. Bufetov, On the Vershik-Kerov Conjecture Concerning the Shannon-McMillan-Breiman Theorem for the Plancherel Family of Measures on the Space of Young Diagrams
Geometric and Functional Analysis. 2012. Vol. 22. No. 4. P. 938-975.

\bibitem{cam-decr} I. Camilier, L. Decreusefond, Quasi-invariance and integration by parts for determinantal and permanental processes,
Journal of Functional Analysis
Volume 259, Issue 1,  2010, pp. 268--300.

\bibitem{DVJ} D.J.Daley, D. Vere-Jones, An introduction to the theory of point processes, vol.I-II, Springer Verlag 2008.

\bibitem{Ghosh}
S. Ghosh,  Determinantal processes and completeness of random exponentials: the critical case, arXiv:1211.2435.

\bibitem{GP}
S. Ghosh, Y. Peres. Rigidity and Tolerance in point processes: Gaussian zeros and Ginibre eigenvalues, arXiv:1211.2381.

\bibitem{G} S. Ghosh,
Rigidity and Tolerance in Gaussian zeros and Ginibre eigenvalues: quantitative estimates, arXiv:1211.3506


\bibitem{HoughEtAl}J.B. Hough, M. Krishnapur, Y. Peres, B. Vir\'ag, Determinantal processes and independence.
Probab. Surv. 3 (2006), 206--229.

\bibitem{IIKS} A.R. Its, A.G. Izergin, V.E. Korepin, N.A. Slavnov, Differential equations for quantum correlation functions, International Journal of Modern Physics B,
Volume 04, Issue 05, April 1990, p.1003-1037.

\bibitem{Kallenberg} O. Kallenberg,  Random Measures. Akademie-Verlag, Berlin,  1986.

\bibitem{Kolmogorov} A. Kolmogoroff, Grundbegriffe der Wahrscheinlichkeitsrechnung, Springer Verlag, 1933.

\bibitem{Lenard} A. Lenard, States of classical statistical mechanical systems of infinitely many particles. I.
Arch. Rational Mech. Anal. 59 (1975), no. 3, 219--239.

\bibitem{Lyons} R. Lyons, Determinantal probability measures.
Publ. Math. Inst. Hautes \'Etudes Sci. No. 98 (2003), 167--212.


\bibitem{Lytvynov} E. Lytvynov, Fermion and boson random point processes as particle distributions of infinite free
Fermi and Bose gases of finite density.
Rev. Math. Phys. 14 (2002), no. 10, 1073--1098.

\bibitem{Macchi} O. Macchi, The coincidence approach to stochastic point processes.
Advances in Appl. Probability, 7 (1975), 83--122.

\bibitem{GO-Adv} G. Olshanski,
The quasi-invariance property for the Gamma kernel determinantal measure. Advances in Mathematics,
2011. Vol. 226. P. 2305-2350

\bibitem{ReedSimon} M. Reed, B. Simon, Methods of modern mathematical physics. vol.I-IV.
Second edition. Academic Press, Inc. New York, 1980.

\bibitem {Rohmes}
Rohlin, V. A. On the fundamental ideas of measure theory. (Russian) Mat. Sbornik N.S. 25(67), (1949),  107 –- 150.

\bibitem{ShirTaka0} T.Shirai, Y. Takahashi, Fermion process and Fredholm
determinant, Proceedings of the Second ISAAC
Congress, vol. I, 15–23, Kluwer 2000.


\bibitem{ShirTaka1}T. Shirai, Y. Takahashi, Random point fields associated with certain Fredholm determinants.
I. Fermion, Poisson and boson point processes. J. Funct. Anal. 205 (2003), no. 2, 414--463.

\bibitem{ShirTaka2}T. Shirai, Y. Takahashi, Random point fields associated with certain Fredholm determinants.
II. Fermion shifts and their ergodic and Gibbs properties. Ann. Probab. 31 (2003), no. 3, 1533--1564.

\bibitem{Simon} B. Simon, Trace class ideals, AMS, 2011.

\bibitem{sinai} Ya.G. Sinai, Theory of phase transitions, rigorous results, Nauka Moscow 1980.

\bibitem{Soshnikov} A. Soshnikov, Determinantal random point fields.
(Russian) Uspekhi Mat. Nauk 55 (2000), no. 5(335), 107--160;
translation in Russian Math. Surveys 55 (2000), no. 5, 923--975.


\bibitem{TW-Airy} C.A. Tracy; H. Widom. Level-spacing distributions and the Airy kernel. Comm. Math. Phys. 159 (1994), no. 1, 151--174.

\bibitem{TracyWidom}C. A. Tracy, H. Widom, Level spacing distributions and the Bessel kernel.
Comm. Math. Phys. 161, no. 2 (1994), 289--309.




\end{thebibliography}
\end{document}